\documentclass[man]{apa6}

\usepackage{pgf,tikz}

\usepackage{apacite}

\usepackage{epsfig}
\usepackage{amssymb}
\usepackage{amsmath}
\usepackage{ifthen}
\usepackage{comment}
\usepackage{setspace}

\title{Counterfactual Graphical Models for Longitudinal Mediation Analysis
with Unobserved Confounding}

\shorttitle{Longitudinal Mediation Analysis}
\rightheader{Longitudinal Mediation Analysis}
\leftheader{I. Shpitser}

\author{
Ilya Shpitser
}
\affiliation{
School of Mathematics\\
University of Southampton\\
\texttt{i.shpitser@soton.ac.uk}
}
\abstract{Questions concerning mediated causal effects are of great interest in
psychology, cognitive science, medicine, social science, public health,
and many other disciplines.  For instance,
about 60\% of recent papers published in leading journals in social psychology
contain at least one mediation test \cite{rucker11mediation}.
Standard parametric approaches to mediation
analysis employ regression models, and either the
``difference method'' \cite{judd81process}, more common in epidemiology,
 or the ``product method'' \cite{baron86mm}, more common in the
social sciences.
In this paper we first discuss a known, but perhaps often unappreciated fact:
that these parametric approaches are a special case of a general counterfactual
framework for reasoning about causality first described by
\citeA{neyman23app},
and \citeA{rubin74potential}, and linked to causal graphical models by
\citeA{robins86new}, and \citeA{pearl00causality}.
We then show a number of
advantages of this framework.  First, it makes the strong
assumptions underlying mediation analysis explicit.  Second, it
avoids a number of problems present in the product and difference methods, such
as biased estimates of effects in certain cases.  Finally, we show
the generality of this framework by proving
a novel result which allows mediation analysis to be applied to longitudinal
settings with unobserved confounders.
}

\def\ci{\perp\!\!\!\perp}
\newcounter{isamac} 
\ifthenelse{\value{isamac}=1}{
\input BoxedEPS          %
\SetTexturesEPSFSpecial  
\HideDisplacementBoxes   %
}{}

\newtheorem{theorem}{Theorem}
\newtheorem{corollary}[theorem]{Corollary}

\newtheorem{definition}[theorem]{Definition}

\newcommand{\path}[1]{{\ensuremath{\boldsymbol{#1}}}} 

\newcommand{\indm}[2]{\ensuremath{{\mathfrak I}_{\kern-1pt\scriptstyle#1}({\mathcal
#2})}} 

\newcommand{\ind}{\mbox{$\perp \kern-5.5pt \perp$}}

\newcommand{\uned}{\hbox{\kern3pt\raise2.5pt\vbox{\hrule
width9pt height 0.3pt}\kern3pt}}

\newcommand{\dashed}{\hbox{\kern3.05pt\raise2.5pt\vbox{\hrule
width1.7pt height 0.3pt}\kern1.8pt\raise2.5pt\vbox{\hrule
width1.7pt height 0.3pt}\kern1.8pt\raise2.5pt\vbox{\hrule
width1.7pt height 0.3pt}\kern1.8pt\raise2.5pt\vbox{\hrule
width1.7pt height 0.3pt}\kern3.05pt}}

\newcommand{\head}{\ensuremath{\succ}}

\newcommand{\pedg}[2]{\ensuremath{{\kern0.5pt
\scriptstyle{\ifthenelse{\equal{\head}{#1}}{\lhead\kern0.5pt}{#1\kern0.5pt}}\joinrel\relbar
\negthinspace\relbar\joinrel{\kern0.5pt #2}\kern0.5pt}}}

\newcommand{\pdots}{\hbox{\kern2.5pt\raise1.5pt\hbox{\ensuremath{\ldots}}\ke
rn2.5pt}}  

\begin{document}
\maketitle
The aim of empirical research in many disciplines is establishing the presence
of effects by means of either randomized trials, or observational
studies if randomization is not possible.
For example, a celebrated success of empirical research in epidemiology is the
discovery of a causal connection between smoking and lung cancer \cite{doll50smoking}.

Once the presence of an effect is established, the precise mechanism of the
effect becomes a topic of interest as well.  A particularly popular type
of mechanism analysis concerns questions of \emph{mediation}, that is to
what extent a given effect of one variable on another is direct, and to
what extent is it mediated by a third variable.
For example, it is known that genetic variants on
chromosome 15q25.1 increase both smoking behavior, and the risk of lung cancer
\cite{tyler12genetic}.  A public health mediation question of interest here is whether
these variants increase lung cancer risk by directly making the patients
susceptible in some way, or whether the risk increase is driven by the increase
in smoking.

In psychology, interest in mediation analysis began partly due to the
influential S-O-R model \cite{woodworth28dynamic}, where causal relationships
between stimulus and response
are mediated by mechanisms internal to an organism, and partly due to the
multi-stage causality present in many theories in psychology (such as
attitude causing intentions, which in turn cause behavior in social
psychology).  Today, mediation questions are ubiquitous in psychology.
Mediation analysis is used to
explicate theories of persuasion \cite{tormala07persuasion}, ease
of retrieval \cite{schwarz91ease}, \cite{tormala07ease}, cognitive priming
\cite{eagly93priming}, developmental psychology \cite{conger90un}, and
explore many other areas.
In fact, about 60\% of recent papers published in leading journals in social
psychology contain at least one mediation test \cite{rucker11mediation}.

A standard approach for mediation analysis involves the use of (linear)
structural equation models, and the so called ``difference method''
\cite{judd81process}, and ``product method'' \cite{baron86mm}.  The first
method, more common
in epidemiology, considers an outcome model both with and without the mediator and takes the difference in the coefficients for the exposure as the measure of the indirect effect.  The second method, more common in the social sciences,
takes as a measure of the indirect effect the product of (i) the coefficient for the exposure in the model for the mediator and (ii) the coefficient for the mediator in the model for the outcome. These methods suffer from a number of problems.  First, interpreting linear regression parameters as causal parameters is not appropriate when non-linearities or interactions are present in the underlying causal mechanism, and can lead to bias \cite{mackinnon93estimating}, \cite{kaufman04further}.  Second, it is not always the case that a regression parameter is interpretable as a causal parameter, even if the parametric structural assumptions of linearity and no interaction hold \cite{robins86new}.  Finally, these methods are not directly applicable to longitudinal settings (where multiple treatments happen over
time) and assume no unmeasured confounding.

The aim of this paper is twofold.  First, we describe recent developments
in the causal inference literature which address the limitations of the
approaches based on linear structural equations
\cite{judd81process}, \cite{baron86mm}.  In particular, we show
that the linear structural equation approach to mediation analysis is a special
case of a more general framework based on potential outcome counterfactuals,
developed by \citeA{neyman23app} and \citeA{rubin74potential},
and extended and linked to non-parametric structural equations and graphical
models by \citeA{robins86new}, and \citeA{pearl00causality}.
We show how this more general framework avoids the difficulties of the linear
structural equations approach, and has additional advantages in making strong
causal assumptions necessary for mediation analysis explicit.  Second,
we use the counterfactual framework to develop novel results which extend
existing mediation analysis techniques to longitudinal settings with some
degree of unmeasured confounding.

Our argument is that to handle increasingly complex mediation questions in
psychology and cognitive science, scientists must necessarily move beyond
the linear structural equation approach, and embrace more general frameworks
for mediation analysis.  The linear structural equation approach is
simply not applicable in complex data analysis settings, and careless
generalizations of this approach will lead to biased conclusions.

The paper is organized as follows.  In section 2, we describe 
mediation analysis based on linear structural equations
in more detail, and describe situations where the use of this method leads to
problems.  In section 3, we introduce a causal inference framework based
on potential outcome counterfactuals and graphical models, and show how this
framework generalizes the linear structural equation framework, and correctly
handles the problems described in section 2.  In section 4, we describe two
motivating examples involving three complications: unobserved confounding,
longitudinal treatments, and path-specific effects, and show how the
counterfactual framework is able to handle these complications with ease.
Section 5 contains the discussion and our conclusions.  The general theory
necessary to solve examples of the type shown in section 4 is
contained in the supplementary materials.

\section{Mediation Analysis Using Linear Structural Equations Models}

The standard mediation setting contains three variables, the cause or treatment
variable, which we will denote by $A$, the effect or outcome variable, which we
will denote by $Y$, and the mediator variable, which we will denote by $M$.
The treatment $A$ is assumed to have an effect on both mediator $M$ and
outcome $Y$, while the mediator $M$ has an effect on the outcome $Y$.
A typical goal of causal inference is establishing the presence of
the \emph{total effect}, or just the causal effect, of $A$ on $Y$.
The goal of mediation analysis is to decompose the total effect into
the \emph{direct effect} of the
treatment $A$ on the outcome $Y$, with the \emph{indirect} or
\emph{mediated effect} of the
treatment $A$ on the outcome $Y$ through the mediator $M$.

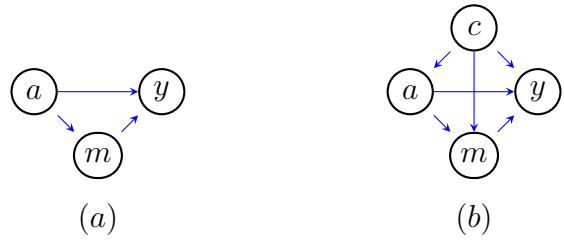
\begin{figure}
\begin{center}
  \begin{tikzpicture}[>=stealth, node distance=1.2cm]

  \begin{scope}
    \tikzstyle{format} = [draw, thick, rectangle, minimum size=6mm, rounded
        corners=3mm]
    \path[->]	node[format] (a) {$a$}
    		node[format, below right of=a] (m) {$m$}
                  (a) edge[blue] (m)
		node[format, above right of=m] (y) {$y$}
		  (m) edge[blue] (y)
		  (a) edge[blue] (y)
		;
    \path[]	node[yshift = -1.7cm, xshift=0.8485cm] (label) {$(a)$};
  \end{scope}
  \begin{scope}[xshift=5.0cm]
    \tikzstyle{format} = [draw, thick, rectangle, minimum size=6mm, rounded
        corners=3mm]
    \path[->]	node[format] (a) {$a$}
    		node[format, below right of=a] (m) {$m$}
                  (a) edge[blue] (m)
		node[format, above right of=m] (y) {$y$}
		  (m) edge[blue] (y)
		  (a) edge[blue] (y)
		node[format, above right of=a] (c) {$c$}
		  (c) edge[blue] (y)
		  (c) edge[blue] (a)
                  (c) edge[blue] (m)
		;
    \path[]	node[yshift = -1.7cm, xshift=0.8485cm] (label) {$(b)$};
  \end{scope}

  \end{tikzpicture}
\end{center}
\caption{(a) Typical mediation setting: $a$ is the treatment, $m$ is
the mediator, $y$ is the outcome.
(b) Mediation setting with observed
confounding: $c$ is a confounder due to being a common cause of treatment,
mediator and outcome variables.
}
\label{triangle}
\end{figure}

Causal relationships in mediation analysis are often displayed by means of
\emph{causal diagrams}.  A causal diagram is a directed graph where nodes
represent variables of interest, in our case the treatment $A$, the mediator
$M$, and the outcome $Y$, and directed arrows represent, loosely, ``direct
causation.''  The mediation setting is typically represented by means of a
causal diagram shown in Fig. \ref{triangle} (a).


The situation represented by this picture contains a treatment that is
either randomly assigned by the experimenter, or randomized naturally.
For example, genetic variants on chromosome 15q25.1 which are linked with
smoking behavior and lung cancer \cite{tyler12genetic} can generally
(modulo possibly some confounding due to population genetics) be assumed to be
naturally randomized.  In psychology, a randomized treatment is often a
treatment or prevention program, such as drug prevention.

Another common situation assumes that the treatment is not randomized, but
all causes of the treatment are observed.  One example of this situation is
shown in Fig. \ref{triangle} (b), which contains a single observed confounder
$C$.  Extending methods described in this section to this case is
straightforward.

Given the causal structure shown in Fig. \ref{triangle}, the statistical
analysis proceeds as follows.  First, the causal relationships between
treatment, mediator and outcome are assumed to take the form of a
causal regression model, or linear structural equation:

\begin{equation}\label{y_model}
Y = \alpha_0 + \alpha_1 \cdot A + \alpha_2 \cdot M + \epsilon_y
\end{equation}
\begin{equation}\label{m_model}
M = \beta_0 + \beta_1 \cdot A + \epsilon_m
\end{equation}
where $\alpha_0, \beta_0$ are intercepts, $\alpha_1,\alpha_2, \beta_1$ are
regression coefficients, $\epsilon_y, \epsilon_m$ are mean zero noise terms,
and the covariance of the noise terms for $Y$ and $M$ is assumed to equal
$0$: $Cov(\epsilon_y, \epsilon_m) = 0$.

For the ``difference method'', a regression model for the outcome where
the mediator is omitted is also included in the analysis:

\begin{equation}\label{y_model_}
Y = \gamma_0 + \gamma_1 \cdot A + \epsilon'_y
\end{equation}

For binary outcomes, it is straightforward to specify alternative regression
models, such as logistic regression models.  However, as we shall soon see,
even this simple modeling change requires care.

The total effect under these models is taken to equal to
$(\alpha_1 + \alpha_2 \cdot \beta_1)$, and can be derived using Sewall Wright's
rules of path analysis \cite{wright21correlation}.
The direct effect under these models is taken to equal to the
regression coefficient $\alpha_1$ of the treatment in the outcome model
(equation \ref{y_model}).
The ``product method'' \cite{baron86mm}, and the ``difference method''
\cite{judd81process} both aim to express the indirect effect of
$A$ on $Y$ in terms of statistical parameters of these regression models.
The product method takes as a measure of the indirect effect the product
of (i) the coefficient for the treatment in the model for the mediator
($\beta_1$ in equation \ref{m_model}), and (ii) the coefficient for the mediator
in the model for the outcome $(\alpha_2$ in equation \ref{y_model}).
The difference method considers the outcome model with (equation \ref{y_model})
and without the mediator (equation \ref{y_model_}), and takes the difference
in the coefficients for the treatment in these two models
($\alpha_1$ and $\gamma_1$) as the measure of the indirect effect.
If the outcome and mediator are continuous and there are no interaction
terms in the regression model for the outcome, the two methods produce
identical answers for the indirect effect \cite{mackinnon95simulation}.

An important property in mediation analysis is the decomposition property:
\begin{equation}\label{decomposition}
\text{Total effect} = \text{Direct Effect} + \text{Indirect Effect}
\end{equation}
This property allows the investigator to quantify how much of an existing
total effect of treatment on outcome is due to the direct influence on
the outcome, and how much is due to the influenced mediated by a third
variable.  Note that it is possible for the total effect to be weak or
non-existent, and direct and indirect effects to both be strong.  This
situation can occur due to cancellation of effects.  For instance, there may
be a strong positive direct effect, but an equally strong negative mediated
effect, resulting in a weak total effect.  The decomposition property holds
for linear structural equation models with continuous outcomes, for indirect
effects defined by both the product and difference methods.

The advantage of the product and difference methods is their simplicity
-- they rely on standard
software for fitting regression models.  The disadvantage is their lack of
flexibility.  In order to work, these methods require assumptions of linearity,
no unmeasured confounding between mediator and outcome, and continuous outcomes.
As we will see in the next section, careless application of
these methods in settings where one or more of these assumptions
are violated will result in bias, and counterintuitive conclusions.


\subsection{Problems with the Product and Difference Methods}

With binary outcomes \cite{mackinnon93estimating}, or interaction terms in the
outcome \cite{tyler09conceptual}, the two methods above no longer
agree on the estimate of the indirect effect.
In addition, there is evidence that in the case of non-linearities or
interactions in the outcome model, neither method gives a satisfactory
measure of the indirect effect \cite{tyler09conceptual}.

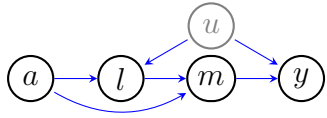
\begin{figure}
\begin{center}
  \begin{tikzpicture}[>=stealth, node distance=1.2cm]

  \begin{scope}
    \tikzstyle{format} = [draw, thick, rectangle, minimum size=6mm, rounded
        corners=3mm]
    \path[->]	node[format] (a) {$a$}
		node[format, right of=a] (l) {$l$}
		  (a) edge[blue] (l)
    		node[format, right of=l] (m) {$m$}
                  (l) edge[blue] (m)
		  (a) edge[blue, bend right] (m)
		node[format, right of=m] (y) {$y$}
		  (m) edge[blue] (y)
		node[format, gray, above of=m, yshift=-0.5cm] (u) {$u$}
		  (u) edge[blue] (l)
		  (u) edge[blue] (y)
		;
  \end{scope}

  \end{tikzpicture}
\end{center}
\caption{Mediation setting with an unobserved confounder $U$, two mediators
$L$ and $M$, and no direct effect of $A$ on $Y$. }
\label{verma}
\end{figure}

Furthermore, even for the case of continuous outcome models with no
interaction terms, certain underlying causal structures can make it impossible
to associate any standard regression parameter with direct and indirect effects.
Consider the causal diagram shown in Fig. \ref{verma}.  This diagram represents
a situation where we have a randomized treatment $A$ and the outcome $Y$,
but instead of a single mediator, we have two mediating variables $L$ and $M$.
Furthermore, we have reasons to believe there is a strong source of unobserved
confounding (which we call $U$) between one of the mediators $L$ and
the outcome $Y$.  For instance, if $A$ represents a primary prevention
program (say drug prevention), and $M$ represents a secondary prevention
program (say a program designed to increase screening rates for serious illness),
then $L$ might represent some observable intermediate outcome of people enrolled
in the primary program, perhaps linked to eventual outcome $Y$ via some
unobserved measure of conscientiousness or health consciousness.
Assume for the moment that all variables are continuous,
and we can model their relationships using linear regression models:

\begin{equation}\label{v_y_model}
Y = \alpha_0 + \alpha_1 \cdot A + \alpha_2 \cdot M + \alpha_3 \cdot L
	+ \epsilon_y
\end{equation}

\begin{equation}\label{v_m_model}
M = \beta_0 + \beta_1 \cdot A + \beta_2 \cdot L + \epsilon_m
\end{equation}

\begin{equation}\label{v_l_model}
L = \delta_0 + \delta_1 \cdot A + \epsilon_l
\end{equation}

We model the presence of the $U$ confounder by allowing that
$Cov(\epsilon_y,\epsilon_l) \neq 0$, while assuming
$Cov(\epsilon_y,\epsilon_m) = 0, Cov(\epsilon_m,\epsilon_l) = 0$.  We still
assume mean zero error terms.
Note that though the directed arrows in the graph in Fig. \ref{verma} are causal,
not all of the regression coefficients in above equations have causal interpretations.
In particular, $\alpha_1$, and $\alpha_3$ do not have causal interpretations, while
$\alpha_2$ does (as the direct effect of $M$ on $Y$).

We are interested in quantifying the direct effect of $A$ on $Y$, and the
effect of $A$ on $Y$ mediated by $M$.  The question is, what (combination of)
parameters of the regression models we specified correspond to these effects.
A naive approach would be to consider a regression model in equation
(\ref{v_y_model}), and take the regression parameter $\alpha_1$
associated with $A$
as the measure of the direct effect.  This approach is wrong, and will lead to
bias.  The difficulty with this example is that a regression coefficient
of a particular independent variable $X$
represents the extent to which the dependent variable $Y$ depends on $X$
given that we condition on all other independent variables.  In our example,
the regression coefficient for $A$ represents dependence of $Y$ on $A$
given that we conditioned on $L$ and $M$ (we do not condition on $U$ since
$U$ is not observed).  Unfortunately, conditioning on $L$ makes $U$ and
$A$ dependent due to the phenomenon known as ``explaining away.''

Consider a toy causal system: a light in a hallway is wired to two
toggle light switches on the opposite ends of the hallway.  If either of
the light switches is flipped, the light turns on.  Two people, Alice and Uma,
stand at opposite ends of the hallway, each near a switch.  Alice sees the
light turn on, and knows she did not toggle the switch.  She can then conclude
(``explain away'' the light turning on)
that Uma toggled the switch.  In our graph, Alice's switch is $A$, Uma's switch
is $U$, and the light itself is $L$.  Conditional on $L$, we can learn
information about $U$ if we know something about $A$.  In other words,
conditional on $L$, $A$ and $U$ become dependent.  Of course, $U$ is a direct
cause of $Y$.  This means that
some of the variation of $Y$ due to $A$, represented by the regression
coefficient of $A$ in equation (\ref{v_y_model}) is actually due to
the ``explaining away'' effect correlating $A$ and $U$, which in turn correlates
$A$ and $Y$ in a non-causal way.  In particular, even if there is no
direct effect of $A$ on $Y$, the regression coefficient of $A$ will not vanish
in most models.


In fact, it can be shown that in examples of this sort,
the presence of unobserved confounders,
coupled with the ``explaining away'' effect will preventing us from associating
any standard function of regression coefficients with causal parameters in a way
that avoids bias.  Furthermore, even if the correct expression for the direct
effect is used (as derived in a subsequent section, and shown in equation \ref{v_dir_eff}),
using standard statistical models in that expression can result in cases where the absence
of direct effect is not possible given the model.  In particular, if we use a linear
regression model with no interaction terms for a continuous outcome $Y$, and a logistic
regression model with no interaction terms for a binary mediator $L$, then the absence of
direct effect is impossible given those models in the sense that the expression in
(\ref{v_dir_eff}) will \emph{never} equal $0$.
This difficulty, which applies not only to regression
models but to almost any standard parametric statistical model associated with
causal diagrams such as Fig. \ref{verma}, is known as the ``null paradox''
\cite{robins97estimation}.

Finally, even if assumptions of linearity, no interaction, and no unobserved
confounding hold, \emph{no} function of the observed data will equal to either
direct or indirect effect in general.  In order
for this equality to hold, it must be the case that error terms of the
outcome and mediator model remain uncorrelated for any possible set of
assignments of independent variables to the model.  This assumption is also
necessary in order to derive mediation effects from double randomization
studies \cite{word74nonverbal} \cite{imai13experimental}.  The assumption
cannot easily be tested, and can be viewed as ruling out unobserved
confounding between variables in different counterfactual situations.
It will be described in more detail later.  Deriving analogues of
this crucial assumption in more complex settings, for the purposes of
sensitivity analysis, can be challenging.

A way out of many of these difficulties involves generalizing from linear
regression models to a general non-parametric framework based on potential
outcome counterfactuals.  This framework will be described in great detail
in the next section.  We will show how this framework gives a more general
representation of direct and indirect effects that will happen to coincide
with the results of the product and difference methods in the special case
of linear regression models.  We will also show how assumptions underlying
mediation analysis can be clearly explained as independence statements
among random counterfactual variables, displayed graphically by a causal
diagram.  We will discuss possible solutions to
the null paradox that can be derived in this framework.
Finally, the flexibility of the framework will allow us to pose more complex
questions of mediation, such as ``what is the effect of $A$ on $Y$ along
the path $A \to M \to Y$ in the graph in Fig. \ref{verma}?'' and answer
these questions in complex settings involving multiple time-dependent
treatments, and unobserved confounding.


\section{Potential Outcomes and Mediation}


Typically, the notion of causal effect of treatment $A$ on outcome $Y$ refers
to change in the outcome between the control group and the test group in
a randomized control trial.  A general representation of causality, divorced
from a particular statistical model such as a regression model must capture this
notion in some way.  An idealized, mathematical representation of a randomized
control trial captures the notion of controlling a variable by means of
an \emph{intervention}.  An intervention on $A$, denoted by $\text{do}(a)$ by
\citeA{pearl00causality}, refers to an operation that fixes the value of $A$
to $a$ regardless of the natural variability of $A$.  An intervention represents
an assignment of treatment to the test group, or a decision to set $A$ to
$a$.  The variation in the outcome
after an intervention is captured by means of an \emph{interventional
distribution}, sometimes denoted by $p(y | \text{do}(a))$.

Crucially, intervening to force $A$ to value $a$ is not the same as observing
that $A$ attains the value $a$, that is: $p(y | a) \neq p(y | \text{do}(a))$.
As an example: ``only Olympic sprinters that can run quickly win gold medals
(observation), therefore I should wear a gold medal to run faster
(intervention)''
\footnote{We want to thank the lesswrong.com community for this example.}.
This is the essence of the common refrain that correlation (statistical
dependence) does not imply causation.

A potential outcome counterfactual refers to the value of a random variable
under a particular intervention $\text{do}(a)$ for a particular unit (individual)
$u$, and is denoted by $Y(a,u)$.  If we wish to average over units in a
particular study, we would obtain a random variable $Y(a)$, representing
variation in the outcome after the intervention $\text{do}(a)$ was performed.
In other words, $Y(a)$ is a random variable with a distribution
$p(y | \text{do}(a))$.

Assume for the moment the simplest mediation setting with variables $A$, $M$,
$Y$, shown in Fig. \ref{triangle}, and assume the causal relationships
between $A$, $M$, and $Y$ can be captured by structural equations shown in
equations (\ref{y_model}), and (\ref{m_model}).  The intervention $\text{do}(a)$
in these systems of equations is represented by replacing the random variable
$A$ in each equation with the intervened value $a$.  Alternatively, if we
augment equations (\ref{y_model}) and (\ref{m_model}) with another equation
for $A$ itself, such as:
\begin{equation}\label{a_model}
A = \epsilon_a
\end{equation}
then the intervention on $A$ can be represented by replacing equation
(\ref{a_model}) by another equation that sets $A$ to a constant $a$.
If interventions are represented in this way, then
the total effect of
$A$ on $Y$, equal to $(\alpha_1 + \alpha_2 \cdot \beta_1)$, can be viewed as
\begin{equation}\label{additivity}
\text{Total Effect} = E[Y(a = 1) - Y(a = 0)] = E[Y(1)] - E[Y(0)]
\end{equation}
In other words, the total effect is the expected difference of outcomes
under two hypothetical interventions.  In one intervention, $A$ is set to $1$,
and in another $A$ is set to $0$.  Note that this definition is non-parametric
in that it does not rely on the model for $Y$ being a linear regression model.
In fact, the definition remains sensible even if we replace the models for $Y$
and $M$ by arbitrary functions:

\begin{equation}\label{y_model_general}
Y = f_y(M,A,\epsilon_y)
\end{equation}
\begin{equation}\label{m_model_general}
M = f_m(A,\epsilon_m)
\end{equation}
\begin{equation}\label{a_model_general}
A = f_a(\epsilon_a)
\end{equation}
These models can be viewed as (non-parametric) structural equations, and
are discussed in great detail by \citeA{pearl00causality}.  The key idea is
that we assume the causal relationship between a variable, say $Y$, and
its direct causes is by means of some unrestricted causal mechanism function
$f_y$.  These structural models can still be modeled by means of causal
diagrams, but are no longer bound by linearity, lack of interactions,
or other parametric assumptions.

\subsection{Direct and Indirect Effects As Potential Outcomes}

Representing direct and indirect effects using potential outcomes is slightly
more involved.  In the case of total effects, the intuition was that $A$ being
set to $0$ represents ``no treatment,'' while $A$ being set to $1$ represents
``treatment,'' and we want to subtract off the expected outcome under no
treatment (the baseline effect) from the expected outcome under treatment.
In the case of direct effects we still would like to subtract off the baseline,
but from the effect that considers only the direct influence of $A$ on $Y$ in
some way.

One approach that preserves the attractive property of decomposition of total
effects into direct and indirect effects proceeds as follows.  We consider
a two stage potential outcome.  In the first stage, we consider for a particular
unit $u$, the value the mediator would take
under baseline treatment $a = 0$: $M(0,u)$.  We then consider
the outcome value of that same unit if the treatment was set to $1$, and mediator
was set to $M(0,u)$: $Y(1,M(0,u),u)$.  In other words, the direct influence of
$A$ on $Y$ for this unit is quantified by the value of the outcome in
a hypothetical situation where we give the
individual the treatment, but also force the mediator variable to behave as if
we did not give the individual treatment.  In graphical terms, this is the
outcome value if active treatment $a = 1$ is only active along the direct
path $A \to Y$, but not active along the path $A \to M \to Y$, since we
force $M$ to behave as if treatment was set to $0$ for the purposes of that
path.  If we average over units, we get
a nested potential outcome random variable: $Y(1,M(0))$.  We define the
direct effect as the difference in expectation between this random variable,
and the baseline outcome:

\begin{equation}\label{direct_effect}
\text{Direct Effect} = E[Y(1,M(0))] - E[Y(0)]
\end{equation}

Note that $E[Y(0)] = E[Y(0,M(0))]$.
The indirect effect is defined similarly, expect we now subtract off the
direct influence of $A$ on $Y$ from the total effect of setting $A$ to $1$:

\begin{equation}\label{indirect_effect}
\text{Indirect Effect} = E[Y(1)] - E[Y(1,M(0))]
\end{equation}
It is not difficult to show that these definition reduce to definitions in
terms of regression coefficients given in the previous sections in the special
case where $Y$ is continuous, $f_y$, and $f_m$ are linear functions with
no interactions, and all $\epsilon$ noise terms are Gaussian.
However, these effect definitions, known as \emph{natural} \cite{pearl01direct} or
\emph{pure} \cite{robins92effects}
are the only sensible definitions of direct and indirect effects currently
known that simultaneously maintain the decomposition property
(\ref{decomposition}),
and apply to arbitrary functions in structural equations (\ref{y_model_general}),
(\ref{m_model_general}), and (\ref{a_model_general}).

\subsection{Assumptions Underlying Mediation Analysis}

  Defining the influence of $A$ on $Y$ for a particular unit $u$ as
$Y(1,M(0,u),u)$ involved a seemingly impossible hypothetical situation,
where the treatment given to $u$ was $0$ for the purposes of the mediator
$M$, and $1$ for the purposes of the outcome $Y$.  In other words, this
situation is a function of multiple, conflicting hypothetical worlds.
In general, no experimental
design is capable of representing this situation unless it is possible to
bring the unit by some means to the pre-intervention state (perhaps by means
of a ``washout period,'' or some other method).  In order to express direct
and indirect effects defined in the previous section as functions of the
observed data, such as regression coefficients, we must be willing to make
certain assumptions that make our impossible hypothetical situation amenable
to statistical analysis.

A typical assumption that makes our situation tractable is expressed in terms
of conditional independence statements on potential outcome counterfactuals:
\begin{equation}\label{assumption}
Y(1,m) \ci M(0)
\end{equation}
where $(X \ci Y)$ stands for ``X is marginally independent of Y'', and
$(X \ci Y | Z)$ stands for ``X is conditionally independent of Y given Z.''

This assumption states that if we happen to have some information on how
the mediator varies after treatment is set to $0$, this does not give us
any information about how the outcome varies if we set the treatment to $1$ and
the mediator to (arbitrary) $m$.  Note that this assumption immediately follows
if we assume independent error terms in a non-parametric structural equation model
defined by (\ref{y_model_general}), (\ref{m_model_general}), (\ref{a_model_general}).
This assumption allows us to perform the following derivation:
\begin{equation}\label{first_product}
p(Y(1,M(0))) = \sum_m p(Y(1,m), M(0) = m) = \sum_m p(Y(1,m)) \cdot p(M(0) = m)
\end{equation}
This derivation expressed our potential outcome as a product of two terms,
with each
of these terms representing variation in a random variable after a well defined
intervention.  This represents progress, since we were able to express a
random variable not typically representable by any experimental design in
terms of results of two well defined randomized trials, one involving $Y$
as the outcome and $A,M$ as treatments, and one involving $M$ as the outcome,
and $A$ as treatment.

Unfortunately, even a single randomized study can be expensive or possibly
illegal to perform on people (if the treatment is harmful), let alone two.
For this reason
a common goal in causal inference is to find ways of expressing interventional
distributions as functions of observed data.  In the causal inference literature
this problem is known as the \emph{identification problem of causal effects}.

As mentioned earlier, the interventional distribution, such as that
corresponding to $Y(a)$, namely $p(y | \text{do}(a))$, is not necessarily
equal to a conditional distribution $p(y | a)$.  Nevertheless, such an
equality holds if there is no unobserved confounders, or common causes
between $A$ and $Y$.  This happens to be the
case in our example.  In terms of potential outcomes, the lack of unobserved
confounding is expressed in terms of the \emph{ignorability assumption}
\[
Y(a) \ci A
\]
In words, this assumption states that if we happen to have information on
the treatment variable, it does not give us any information about the outcome
$Y$ after the intervention $\text{do}(a)$ was performed.  A graphical way of
describing ignorability is to say that there does not exist certain kinds of
paths between $A$ and $Y$, called \emph{back-door paths}
\cite{pearl00causality},
in the causal diagram.  Such paths are called ``back-door'' because they start
with an arrow pointing into $A$.
It can be shown that if ignorability holds for $Y(a)$ and $A$ (alternatively
if there are no back-door paths from $A$ to $Y$ in the corresponding causal
diagram), then $p(y | \text{do}(a)) = p(y | a)$.

If there exist common causes of $A$ and $Y$ but they are observed, as is the
case of node $C$ in Fig. \ref{triangle} (b), it is possible to express a
more general assumption known as the \emph{conditional ignorability assumption}
\[
Y(a) \ci A | C
\]
In words, this assumption states that if we happen to have information on
the treatment variable, then conditional on the observed confounder $C$,
this information gives us no information about the outcome $Y$ after the
intervention $\text{do}(a)$ was performed.  In graphical terms, this assumption
is equivalent to stating that $C$ ``blocks'' all back-door paths from $A$ to $Y$
\footnote{\citeA{pearl88probabilistic,pearl00causality} gives
a more detailed discussion of the
notion of ``blocking'' that has to be employed here.}.
It can be shown that if conditional ignorability $(Y(a) \ci A | C)$ holds, then
$p(y | \text{do}(a)) = \sum_c p(y | a,c) p(c)$.  This formula is known as the
back-door formula, or the adjustment formula.

Sometimes,
identification of causal effects is possible even in the presence of unobserved
confounding.  See the work of \citeA{tian02general}, \citeA{huang06do}, and
\citeA{shpitser06id, shpitser06idc,shpitser07hierarchy}
for a general treatment of the causal effect identification problem.

In our case, the ignorability assumption for $M(a)$ and $A$, as well as for
$Y(1,m)$ and $A,M$ allows us to further
express each of the terms in the product in (\ref{first_product}) in
terms of observed data as follows:
\[
\sum_m p(Y(1,m)) \cdot p(M(0) = m) = \sum_m p(Y | A = 1, m) \cdot p(m | A = 0)
\]
Plugging this last expression into the formula (\ref{direct_effect})
for direct effects gives us
\[
\sum_m \left\{ E[Y | A = 1, m] - E[Y | A = 0, m] \right\} p(m | A = 0)
\]
This expression is known as the \emph{mediation formula} \cite{pearl11cmf}.
Note that the mediation formula does not require a particular functional
form for causal mechanisms relating $Y$, $M$ and $A$.

Note also that assumption (\ref{assumption}) is untestable, since it is
positing a marginal independence between two potential outcomes, one of which
involves the treatment being set to $1$, and another involves the
treatment being set to $0$.  A form of this assumption is still necessary in
order to equate direct and indirect effects with functions of regression
coefficients in the simple linear regression setting described in previous
sections, in the sense that violations of the assumption
will generally prevent us from uniquely expressing a given
direct or indirect effect as a function of observed data (e.g. the effect
becomes \emph{non-identifiable}.)  
For this reason, even in the simplest mediation problems, care must
be taken to either justify assumption (\ref{assumption}) on strong substantive
grounds, perform a reasonable sensitivity analysis \cite{tchetgen12semi},
or reduce the mediation problem to a testable problem involving interventions
without conflicts \cite{robins10alternative}.

\subsection{Mediation with Unobserved Confounding}

One of the advantages of the potential outcome framework is its flexibility.
Since it does not rely on parametric assumptions, it can be readily extended
to handle modeling complications.  Consider again our two mediator example
in Fig. \ref{verma}.  We mentioned in the previous section that
product and difference methods will result in biased estimates of direct
effects of $A$ on $Y$ not through $M$, due to a combination
of unobserved confounding and the ``explaining away'' effect in that
example.  A non-parametric
definition of direct effect based on potential outcomes avoids these
difficulties.  Our expression for direct effect is
$E[Y(1,M(0))] - E[Y(0)]$.  Since $A$ is randomized (there is no unobserved
confounding between $A$ and the outcome $Y$), the second term can be shown to
equal $E[Y | A = 0]$.  The first term can be shown, given assumption
(\ref{assumption}), and a general theory of identification of causal
effects \cite{tian02general},\cite{shpitser06id,shpitser07hierarchy}
to equal
\[
E[Y(1,M(0))] = \sum_m \left( \sum_l E[Y|m,l,A=1]p(l|A=1) \right) p(m | A=0)
\]
The direct effect is then equal to
\begin{equation}\label{v_dir_eff}
\text{Direct Effect} = \sum_m
	\left( \sum_l E[Y|m,l,A=1]p(l|A=1) \right) p(m | A=0) - E[Y | A = 0]
\end{equation}
while the indirect effect is equal to
\[
\text{Indirect Effect} = E[Y | A = 1] -
	\sum_m \left( \sum_l E[Y|m,l,A=1]p(l|A=1) \right) p(m | A=0)
\]
It can be shown that not only do the direct and indirect effects
add up to the total effect in this case, but the quantity
(\ref{v_dir_eff}) equals $0$ precisely when the effect of $A$
on $Y$ along the arrow $A \to Y$ is in some sense absent.\footnote{As long
as the parametric models for the functionals in the formula are general enough
to avoid the ``null paradox'' issue.  Linear regressions for all terms suffice,
but a no-interaction linear regression for a continuous outcome $Y$,
and a no-interaction logistic regression for a binary mediator $L$ does not
suffice.  The general rule of thumb is the models must be general enough to
permit the above mean differences to equal zero for some parameter settings.}
However, even though we used simple linear regression models in this example,
neither of these effects reduces to any straightforward function of the
regression coefficients.
It is possible to express these kinds of functional as functions of regression
coefficients in an appropriately adjusted model
(such as the marginal structural model, which is estimated by fitting weighted
regression models \cite{robins00marginal}), or as functions of parameters
in a non-standard parameterization of causal models, where statistical parameters
correspond to causal parameters directly \cite{shpitser11eid},\cite{richardson12nested}.

\section{Path-specific Longitudinal Mediation with Unobserved Confounding}

In the previous section we saw how the presence of unobserved confounders
and multiple mediators can easily result in situations where regression
coefficients cannot be meaningfully associated with direct and indirect causal
effects.  In this section, we consider even more complex mediation settings,
which can nevertheless be handled appropriately using the potential outcome
counterfactual framework representing (possibly non-linear) structural equations.
We motivate the discussion with two examples, one from HIV research, and one
from psychology.

The human immunodeficiency virus (HIV) causes AIDS by attacking and destroying
helper T cells.  If the concentration of these cells falls below a critical
threshold, cell-mediated immunity is lost, and the patient eventually dies to
an opportunistic infection.  Patients infected with HIV with reduced
T cell counts are typically put on courses of anti retroviral therapy (ART),
as a first line therapy.  Unfortunately, side effects of many types of ART
medication may cause poor adherence to the therapy (that is, patients do not
always take the medication on time, or stop taking it altogether).  Side effects
are often caused by toxicity of the medication, or patient's adverse reaction
to the medication.  Severity of the side effects is often linked to the
patient's ``overall health level'' (ill defined, and thus not measured),
which also affects the eventual outcome of
the therapy (survival or death).  If the ART happens to not be very effective
at viral suppression, and results in patient deaths, this could be because the
ART itself is not very good, or it could be due to poor patient adherence.
In other words, poor outcomes of ART results in a natural mediation question
in HIV research -- is the poor total effect possibly
due to cancellation of a strong direct effect of the medication on survival
by an equally strong indirect effect of poor adherence?

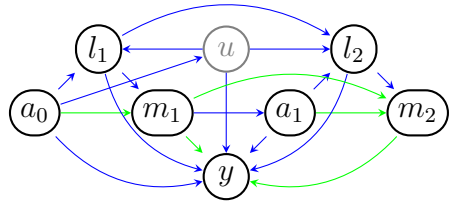
\begin{figure}
\begin{center}
  \begin{tikzpicture}[>=stealth, node distance=1.2cm]

  \begin{scope}
    \tikzstyle{format} = [draw, thick, rectangle, minimum size=6mm, rounded
        corners=3mm]
    \path[->]	node[format] (a0) {$a_0$}
		node[format, above right of=a0] (l1) {$l_1$}
		  (a0) edge[blue] (l1)
    		node[format, below right of=l1] (m1) {$m_1$}
                  (l1) edge[blue] (m1)
		  (a0) edge[green] (m1)
		node[format, gray, above right of=m1] (u) {$u$}
		  (u) edge[blue] (l1)
		node[format, below right of=u] (a1) {$a_1$}
		  (m1) edge[blue] (a1)
		node[format, above right of=a1] (l2) {$l_2$}
		  (a1) edge[blue] (l2)
    		node[format, below right of=l2] (m2) {$m_2$}
		  (l2) edge[blue] (m2)
		  (a1) edge[green] (m2)
		node[format, below right of=m1] (y) {$y$}
		  (m2) edge[green, bend left] (y)
		  (m1) edge[green] (y)
		  (a1) edge[blue] (y)
		  (a0) edge[blue, bend right] (y)
		  (l1) edge[blue, bend right] (y)
		  (l2) edge[blue, bend left] (y)

		  (a0) edge[blue] (u)
		  (u) edge[blue] (l2)
		  (u) edge[blue] (y)

		  (l1) edge[blue, bend left] (l2)
		  (m1) edge[green, bend left=25] (m2)

		;
  \end{scope}

  \end{tikzpicture}
\end{center}
\caption{
A longitudinal mediation setting with an unobserved confounder $u$,
where we are interested in effects along
all paths from $A_0$ and $A_1$ to $Y$ only through $M_1$ or $M_2$.
Paths we are interested in are shown in green.
}
\label{time_med}
\end{figure}

The situation is shown graphically in Fig. \ref{time_med}.  Here, we show
ART taken over the course of two months, represented by two time slices.
In practical studies, ART is taken over a period of years, and the number of
time slices is quite large.  In this graph, the ART is represented by
nodes $A_0$ and $A_1$, the patient outcome by $Y$, patient adherence at each
time slice by $M_1$ and $M_2$, toxicity of the medication by $L_1$, $L_2$,
and finally the unobserved state of patient's health affecting reaction to
the medication and the outcome by $U$.  Since we are interested in the indirect
effect of ART on survival mediated by adherence, we are only interested in
the effect along the paths from $A_0,A_1$ on $Y$ which pass through $M_1,M_2$.
These paths are shown in green in the graph.

Our second example, which is isomorphic to the HIV example above,
concerns the use of prevention programs, to promote
positive outcomes in vulnerable populations.  Assume the primary interventions
$A_0,A_1$ involve attending a drug prevention program.  Like ART in the previous
example, the program is an ongoing intervention (say a monthly meeting).
However, there is also a secondary intervention $M_1,M_2$ which is meant to
increase the rates of screening for serious illness such as cancer.  Although
the secondary intervention is not directly related to drug prevention,
it is conceivable that there is
a synergistic effect between the primary and secondary intervention to promote
positive outcomes (say staying drug free), perhaps due to the fact that
both interventions promote good habits and health consciousness.  In this
example, $L_1,L_2$ are, loosely speaking, ``the participant's responsiveness''
which may be affected by unobserved factors involving family, friends,
socioeconomic background ($U$), and so on.  These unobserved factors also
influence the outcome.  The mediation question here is quantifying the extent
to which the outcome is influenced by the primary intervention itself, versus
an indirect effect via the secondary intervention.  The indirect effect mediated
by the secondary intervention is, again, shown in green in Fig. \ref{time_med}.

Aside from unobserved confounding represented by $U$, a complication also
present in the example in Fig. \ref{verma} in the previous section,
what makes these
examples difficult is first the longitudinal setting where treatments recur
over multiple time slices, and second that we are interested in effects along
a particular bundle of causal paths.  In previous sections we were interested
in effects either only along the direct path $A \to Y$, or only along all
paths other than the direct one.  In these cases we are interested in some
indirect paths (through $M_1, M_2$), but not others (through $L_1, L_2$).

In the case of multiple treatments, causal effect from these treatments
to the outcome $Y$ is transmitted along special paths called
\emph{proper causal paths} of the form
$A_k \to \ldots \to Y$, where $A_k$ is one of the treatments, and where
this path cannot intersect any other treatment other than $A_k$
(otherwise it is really a causal
path from that second treatment to the outcome).  We are interested in
quantifying the effect along a subset of such paths, displayed graphically
as those paths consisting entirely of green arrows.  Algebraically, we will
denote this set of paths of interest as $\pi$.  The proper causal
paths along which
the causal effect is transmitted but which do not lie in $\pi$ is displayed
graphically as those paths which contain at least one blue arrow.

\subsection{Formalization of Path-Specific Effects}

Naturally, even if the statistical model associated with the causal diagram
shown in Fig. \ref{time_med} is given in terms of linear regressions, it is not
possible to express effects of interest as simple functions of regression
coefficients.  However, it is possible to express these \emph{path-specific
effects} \cite{pearl01direct},\cite{chen05ijcai} in terms of potential outcome counterfactuals.

We will use an inductive rule to construct a potential outcome representing
the effect of $A_0,A_1$ on $Y$ only along green paths.  For the purposes of
this rule, we will represent values
$a'_0,a'_1$ of $A_0,A_1$ to represent ``baseline treatment'' or ''no
treatment'' (in the previous section we used $0$),
and values $a_0,a_1$ to represent ``active treatment'' (in the previous
section we used $1$).
This potential outcome will involve $Y$, and interventions on all direct
causes of $Y$ (that is all nodes $X$ such that $X \to Y$ exists in the graph).
These interventions are defined as follows.

If the arrow $X \to Y$ is blue, this means we are not interested in the effect
transmitted along this arrow.  In previous sections we represented this by
considering the value of $X$ ``as if treatment was baseline,'' or $X(0)$.
In our case we will do the same, except since we have two treatments
we set them both to the baseline values $a'_0,a'_1$.
For example, $L_1$ and $L_2$ are both direct causes
of $Y$ along blue arrows.  This means we intervene on whatever values they would
have had if $A_0,A_1$ were set to baseline, or $L_1(a'_0,a'_1)$,
$L_2(a'_0,a'_1)$.
If the direct cause of $Y$
is one of the treatments $A_0$ or $A_1$ then we intervene
to set their value to ``active'' if the arrow from treatment to outcome is
green (e.g.  we are interested in the effect), or ``baseline'' if the arrow
is blue (e.g. we are not interested in the effect).  Finally, if the direct
cause of $Y$ is not a treatment, but is a direct cause along a green arrow,
we inductively set the value of that cause to whatever value
it would have had under a path-specific effect of $A_0,A_1$ on that cause.
For instance, $M_1$ is a direct cause of $Y$ along a green arrow, which means
we set the value of $M_1$ to whatever value is dictated by the path-specific
effect of $A_0,A_1$ on $M_1$.  To figure out what that value is, we simply
apply our rule inductively from the beginning, except to $M_1$ as the outcome,
rather than $Y$.

Applying the first stage of our rule gives us the potential outcome\\
$Y(a'_0, a'_1, L_1(a'_0), L_2(a'_0,a'_1), \gamma_{m_1}, \gamma_{m_2})$,
where $\gamma_{m_1}$ and $\gamma_{m_2}$ are path-specific effects of
$A_0,A_1$ on $M_1$ and $M_2$, respectively, along the green paths only.
If we apply the rule inductively $\gamma_{m_1}$, and $\gamma_{m_2}$, and plug
in and simplify, we get that our path-specific effect is the following complex
potential outcome
\begin{equation}\label{path_spec}
Y(a'_0, a'_1, L_1(a'_0), L_2(a'_0,a'_1), M_1(a_0,L_1(a'_0)),
M_2(a_1,L_2(a'_1,a'_0)))
\end{equation}
While this expression looks algebraically complex, what it is actually
expressing is
a rather simple idea.  We have two treatment levels: ``baseline'' and ``active.''
For the purposes of green paths, the causal paths we are interested in, we
pretend treatment levels are active.  For the purposes of all other paths,
we pretend treatment levels are baseline.  In this way, the treatment is
active only along the paths we are interested in, and all other paths are
``turned off.''  We use this rule to select what values we intervene on, and
then use these interventions in a nested way, following the causal paths
of the graph.  An equivalent definition of path-specific effects
phrased in terms of replacing structural
equations is given by \citeA{pearl01direct}.

\subsection{The Total Effect Decomposition Property for Path-Specific Effects}

In the previous sections we defined the direct and indirect effects by taking
a difference of expectations (see equations \ref{direct_effect} and
\ref{indirect_effect}).  We can generalize such definitions to path-specific
effects to obtain a decomposition of the total effect into a sum of two terms,
one representing the effect along proper causal paths in $\pi$, and one
representing the effect along proper causal paths not in $\pi$.

First, assume the distribution for the nested potential outcome
defining the path-specific effect along proper causal paths in $\pi$
of $A$ on $Y$ is given by $p_{\pi}(Y)$.  Then we have

\begin{equation}\label{path_specific_pi}
\text{Effect along paths in $\pi$} =
	E[Y]_{p_{\pi}(Y)} - E[Y(a')]
\end{equation}
and
\begin{equation}\label{path_specific_not_pi}
\text{Effect along paths not in $\pi$} =
	E[Y(a)] - E[Y]_{p_{\pi}(Y)}
\end{equation}
Since the total effect can be defined as $E[Y(a)] - E[Y(a')]$, we have
\[
\text{Total Effect} = \text{Effect along paths in $\pi$} + \text{Effect along
paths not in $\pi$}
\]

which is an intuitive additivity property stating that the total effect
can be decomposed into a sum of two terms, where one term quantifies the
effect operating along a given bundle of proper causal paths $\pi$, and another
term quantifies the effect operating along all proper causal paths other than those in $\pi$.
Note that this property generalizes the additivity property for direct and
indirect effects, where $\pi$ was taken to mean a single arrow from $A$ to $Y$.

\subsection{Path-Specific Effects as Functions of Observed Data}

In the previous section we showed that a path-specific effect can be
defined in terms of a nested potential outcome after we ``turn off''
causal paths we are not interested in.  Regardless of how sensible such
a definition may be, it is not very useful unless this potential outcome
can be expressed as a function of the observed data, and thus become amenable
to statistical analysis.

In a previous section, we showed that in order to express direct and indirect
effects in terms of observed data, we needed to make an untestable independence
assumption (shown in equation \ref{assumption}).  Path-specific effects
generalize direct and indirect effects, and thus require even more assumptions.

In fact, for the path-specific effect along green paths in Fig. \ref{time_med},
it suffices to believe the following independence claim for any value assignments
$a_0,a_1,a'_0,a'_1,l_1,l_2$:
\begin{equation}\label{cross_world}
\left\{ Y(a'_0,a'_1,l_1,l_2,m_1,m_2), L_1(a'_0), L_2(a'_0,a'_1) \right\} \ci
\left\{ M_1(a_0,l_1), M_2(a_0,a_1,l_2) \right\}
\end{equation}
If we believe this assumption, we can express the path-specific effect in
equation \ref{path_spec} as
\begin{equation}
\label{pot_out_example}
\sum_{l_1,l_2,m_1,m_2} 
	p(Y(a'_0,a'_1,l_1,l_2,m_1,m_2), L_1(a'_0) = l_1,
	L_2(a'_0,a'_1) = l_2) \cdot
\end{equation}
\begin{equation*}
	p(M_1(a_0,l_1) = m_1, M_2(a_1,l_2) = m_2)
\end{equation*}

This expression is a product of terms, where each term is a well defined
interventional density (that is, there are no conflicts involving different
hypothetical worlds).
If we further make use of the general theory of identification
of interventional densities from observed data \cite{tian02on},
\cite{shpitser06id,shpitser07hierarchy},
we can express the above as follows

\begin{equation} \label{observational_g}
\sum_{l_1,l_2,m_1,m_2}
	p(Y | a'_0,a'_1,l_1,l_2,m_1,m_2) \cdot
	p(m_2 | l_2, a_1, m_1, a_0) \cdot
	p(l_2 | a'_0, a'_1, l_1) \cdot
	p(m_1 | l_1, a_0) \cdot
	p(l_1 | a'_0)
\end{equation}
This expression is a function of the observed data.

What is left is finding an expression for the total effect of
$A_0,A_1$ on $Y$.  It is not difficult to show that
$E[Y(a_0,a_1)] = \sum_{m_1,l_1} E[Y | m_1,l_1,a_1,a_0] p(m_1,l_1 | a_0)$.
By analogy with a previous section, we can express the path-specific effect
via fully green paths in $\pi$ through $M_1, M_2$ as a difference of expectations
\begin{align}\label{dir_diff_gen}
\sum_{l_1,l_2,m_1,m_2}
\notag
	E[Y | a'_0,a'_1,l_1,l_2,m_1,m_2] \cdot
	p(m_2 | l_2, a_1, m_1, a_0) \cdot
	p(l_2 | a'_0, a'_1, l_1) \cdot
	p(m_1 | l_1, a_0) \cdot
	p(l_1 | a'_0) -\\
\sum_{m_1,l_1} E[Y | m_1,l_1,a'_1,a'_0] p(m_1,l_1 | a'_0)
\end{align}
while the effect via all paths that are not fully green (that is proper causal
paths not in $\pi$) as another difference
\begin{align}\label{indir_diff_gen}
\notag
\sum_{m_1,l_1} E[Y | m_1,l_1,a_1,a_0] p(m_1,l_1 | a_0) -\\
\sum_{l_1,l_2,m_1,m_2}
	E[Y | a'_0,a'_1,l_1,l_2,m_1,m_2] \cdot
	p(m_2 | l_2, a_1, m_1, a_0) \cdot
	p(l_2 | a'_0, a'_1, l_1) \cdot
	p(m_1 | l_1, a_0) \cdot
	p(l_1 | a'_0)
\end{align}
These quantities can be estimated with standard statistical methods by
simply positing a model for each term, for instance a regression model,
estimating the models from data, and computing the estimated functional.
This is the so called ``parametric g-formula'' approach \cite{robins86new}.
With this method, care must be taken to avoid the ``null paradox'' issue,
as was the case with direct and indirect effects.
A less straightforward approach which only relies on modeling the probability
of the treatment in each times lice given the past,
is to generalize marginal structural models \cite{robins00marginal}, which were
originally developed in the context of estimating total effects in longitudinal
settings with confounding.  Another alternative is to extend existing multiply
robust methods for point treatment mediation settings \cite{tchetgen12semi2} based on
semi-parametric statistics to the longitudinal setting.

\subsection{Expressing Arbitrary Path-Specific Effects in Terms of
	Observed Data}

One difficulty with path-specific effects is that the corresponding
potential outcome counterfactual is nested, and therefore complicated.
On the other hand, the graphical representation of path-specific effects on
a causal diagram is fairly intuitive (effect along green paths only).
For this reason, it would be desirable to obtain a result which says, for
a particular bundle of green paths on a particular causal diagram, whether
the corresponding counterfactual can be expressed as a function of the observed
data, without going into the details of the counterfactual itself.  In this
section we give just such a result, which generalizes existing
results on path-specific effects in cases with a single treatment and no
unobserved confounding \cite{chen05ijcai}.

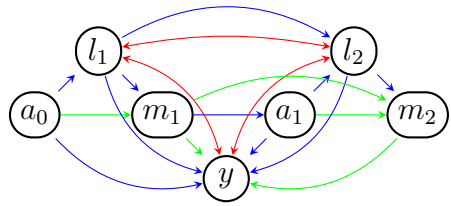
\begin{figure}
\begin{center}
  \begin{tikzpicture}[>=stealth, node distance=1.2cm]

  \begin{scope}
    \tikzstyle{format} = [draw, thick, rectangle, minimum size=6mm, rounded
        corners=3mm]
    \path[->]	node[format] (a0) {$a_0$}
		node[format, above right of=a0] (l1) {$l_1$}
		  (a0) edge[blue] (l1)
    		node[format, below right of=l1] (m1) {$m_1$}
                  (l1) edge[blue] (m1)
		  (a0) edge[green] (m1)
		node[format, below right of=m1] (y) {$y$}
		node[format, above right of=y] (a1) {$a_1$}
		  (m1) edge[blue] (a1)
		node[format, above right of=a1] (l2) {$l_2$}
		  (a1) edge[blue] (l2)
    		node[format, below right of=l2] (m2) {$m_2$}
		  (l2) edge[blue] (m2)
		  (a1) edge[green] (m2)
		  (m2) edge[green, bend left] (y)
		  (m1) edge[green] (y)
		  (a1) edge[blue] (y)
		  (a0) edge[blue, bend right] (y)
		  (l1) edge[blue, bend right] (y)
		  (l2) edge[blue, bend left] (y)

		  (l1) edge[<->,red, bend left=10] (l2)
		  (l1) edge[<->,red, bend left] (y)
		  (l2) edge[<->,red, bend right] (y)

		  (l1) edge[blue, bend left] (l2)
		  (m1) edge[green, bend left=25] (m2)

		;
  \end{scope}
  \end{tikzpicture}
\end{center}
\caption{
A mixed graph representing unobserved confounding in Fig. \ref{time_med}.
}
\label{mixed_time_med}
\end{figure}

We first start with a few preliminaries on graphs.
We will display causal diagrams with unobserved confounding, such that those
in Figs. \ref{verma} and \ref{time_med} by means of a special kind of mixed
graph containing two kinds of edges, directed edges (${\to}$), either blue or
green depending on whether we are interested in the corresponding causal path,
and red bidirected edges (${\color{red}\leftrightarrow}$).  The former represent
direct causation edges, as before.  The latter represent the presence of some
unspecified unobserved common cause.  For example, we represent the causal
diagram in Fig. \ref{time_med} by means of the mixed graph shown in Fig.
\ref{mixed_time_med}.  Note that since $U$ links three nodes, $L_1, L_2, Y$,
each pair of these three is joined by a bidirected arrow.  We call this type
of mixed graph an acyclic directed mixed graph (ADMG) \cite{richardson09factorization}.
\citeA{verma90equiv} called these types of graphs latent projections.
The reason we use bidirected arrows is both to avoid cluttering the graph
with potentially many possible unobserved confounders, and because certain
crucial definitions involving confounders are easier to state in terms of
bidirected arrows.

A sequence of distinct edges such that the first edge connects to $X$, the
last edge connects to $Y$, and each $k$th and $k+1$th edge shares a single
node in common is called a \emph{path} from $X$ to $Y$.
If vertices $X,Y$ are
connected by a path of the form $X \to \ldots \to Y$, then we say $X$ is an
ancestor of $Y$ and $Y$ is a descendant of $X$.  Such a path is called
directed.  A path is called bidirected if it consists exclusively of
bidirected edges.  For an ADMG ${\cal G}$ with a set of nodes $V$, we define
a subgraph ${\cal G}_A$ over a subset $A \subseteq V$ of nodes to consist only
of nodes in $A$ and edges in ${\cal G}$ with both endpoints in $A$.
If an edge $X \to Y$ exists, $X$ is called a parent of $Y$, and $Y$ a child of
$X$.

\begin{definition}[district]
Let ${\cal G}$ be an ADMG.
Then for any node $a$, the set of nodes in ${\cal G}$
reachable from $a$ by bidirected paths is called the
district of $a$, written $Dis_{\cal G}(a)$.
\end{definition}

For example, in the graph in Fig. \ref{mixed_time_med}, $Dis_{\cal G}(Y) =
\{ Y, L_1, L_2 \}$.

\begin{definition}[recanting district]
Let ${\cal G}$ be an ADMG, $A$, $Y$ sets of nodes in ${\cal G}$, and
$\pi$ a subset of proper causal paths which start with a node in $A$ and end in
a node in $Y$ in ${\cal G}$.
Let $V^*$ be the set of nodes not in $A$ which are ancestral of $Y$
via a directed path which does not intersect $A$.  Then a district $D$ in
an ADMG ${\cal G}_{V^*}$ is called
a \emph{recanting district} for the $\pi$-specific effect of $A$ on $Y$
if there exist nodes $z_i, z_j \in D$ (possibly $z_i = z_j$), $a_i \in A$,
and $y_i,y_j \in Y$ (possibly $y_i = y_j$) such that there is
a proper causal path $a_i \to z_i \to \ldots \to y_i$ in $\pi$, and a proper
causal path $a_i \to z_j \to \ldots \to y_j$ not in $\pi$.
\end{definition}

It turns out that the recanting district criterion characterizes situations
when a potential outcome counterfactual can be expressed in terms of
well-defined interventions, without conflicts, as long as we assume that
the causal diagram represents a set of (arbitrary) structural equations.

\begin{theorem}[recanting district criterion]
\label{full-result}
Let ${\cal G}$ be an ADMG representing a causal diagram with unobserved
confounders corresponding to a structural causal model.
Let $A$, $Y$ sets of nodes nodes in ${\cal G}$, and
$\pi$ a subset of proper causal paths which start with a node in $A$ and end
in a node in $Y$ in ${\cal G}$.
Then the $\pi$-specific effect of $A$ on $Y$ is expressible as a functional of
interventional densities if and only if there does
not exists a recanting district for this effect.
\end{theorem}

The functional referenced in the theorem is equal to
\begin{equation}\label{inter_g}
\sum_{V^* \setminus Y} \prod_{D} p(D = d | \text{do}(E_D = e_D))
\end{equation}
where D ranges over all districts in the graph ${\cal G}_{V^*}$,
$E_D$ refers to nodes with directed arrows pointing into $D$ but which are
themselves not
in $D$, and value assignments $d,e_D$ are assigned as follows.  If any element
$a$ in $A$ occurs in $E_D$ in a term $p(D = d | \text{do}(E_D = e_D))$, then
it is assigned a baseline value if the arrows from $a$ to elements in $D$ are
all blue, and an active value if the arrows from $a$ to elements in $D$ are all
green.  All other elements in $E_D$ are assigned values consistent with the
values indexed in the summation.

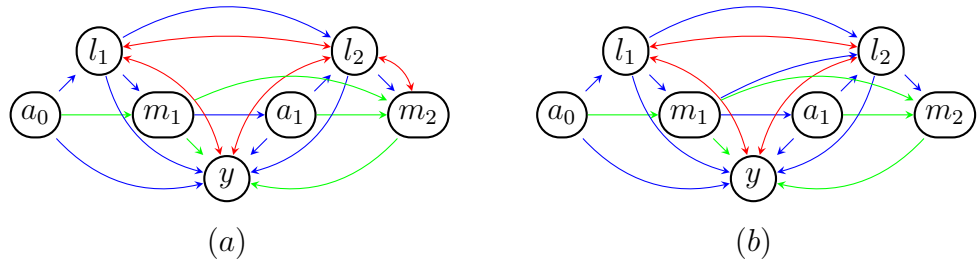
\begin{figure}
\begin{center}
  \begin{tikzpicture}[>=stealth, node distance=1.2cm]

  \begin{scope}
    \tikzstyle{format} = [draw, thick, rectangle, minimum size=6mm, rounded
        corners=3mm]
    \path[->]	node[format] (a0) {$a_0$}
		node[format, above right of=a0] (l1) {$l_1$}
		  (a0) edge[blue] (l1)
    		node[format, below right of=l1] (m1) {$m_1$}
                  (l1) edge[blue] (m1)
		  (a0) edge[green] (m1)
		node[format, below right of=m1] (y) {$y$}
		node[format, above right of=y] (a1) {$a_1$}
		  (m1) edge[blue] (a1)
		node[format, above right of=a1] (l2) {$l_2$}
		  (a1) edge[blue] (l2)
    		node[format, below right of=l2] (m2) {$m_2$}
		  (l2) edge[blue] (m2)
		  (a1) edge[green] (m2)
		  (m2) edge[green, bend left] (y)
		  (m1) edge[green] (y)
		  (a1) edge[blue] (y)
		  (a0) edge[blue, bend right] (y)
		  (l1) edge[blue, bend right] (y)
		  (l2) edge[blue, bend left] (y)

		  (l1) edge[<->,red, bend left=10] (l2)
		  (l1) edge[<->,red, bend left] (y)
		  (l2) edge[<->,red, bend right] (y)
		  (l2) edge[<->, red, bend left] (m2)

		  (l1) edge[blue, bend left] (l2)
		  (m1) edge[green, bend left=25] (m2)

		;
    \path[]	node[yshift = -1.7cm, xshift=2.5485cm] (label) {$(a)$};
  \end{scope}
  \begin{scope}[xshift=7.0cm]
    \tikzstyle{format} = [draw, thick, rectangle, minimum size=6mm, rounded
        corners=3mm]
    \path[->]	node[format] (a0) {$a_0$}
		node[format, above right of=a0] (l1) {$l_1$}
		  (a0) edge[blue] (l1)
    		node[format, below right of=l1] (m1) {$m_1$}
                  (l1) edge[blue] (m1)
		  (a0) edge[green] (m1)
		node[format, below right of=m1] (y) {$y$}
		node[format, above right of=y] (a1) {$a_1$}
		  (m1) edge[blue] (a1)
		node[format, above right of=a1] (l2) {$l_2$}
		  (a1) edge[blue] (l2)
    		node[format, below right of=l2] (m2) {$m_2$}
		  (l2) edge[blue] (m2)
		  (a1) edge[green] (m2)
		  (m2) edge[green, bend left] (y)
		  (m1) edge[green] (y)
		  (a1) edge[blue] (y)
		  (a0) edge[blue, bend right] (y)
		  (l1) edge[blue, bend right] (y)
		  (l2) edge[blue, bend left] (y)

		  (l1) edge[<->,red, bend left=10] (l2)
		  (l1) edge[<->,red, bend left] (y)
		  (l2) edge[<->,red, bend right] (y)
		  (m1) edge[blue, bend left=10] (l2)

		  (l1) edge[blue, bend left] (l2)
		  (m1) edge[green, bend left=25] (m2)

		;
    \path[]	node[yshift = -1.7cm, xshift=2.5485cm] (label) {$(b)$};
  \end{scope}

  \end{tikzpicture}
\end{center}
\caption{ Two variations on the graph in Fig. \ref{mixed_time_med} where
the path-specific effect of $A_0,A_1$ on $Y$ along the green paths is not
identifiable from observational (or even interventional) data.
(a) The presence of an unobserved parent of $L_2$ and $M_2$ spoils
identification. (b) If $M_1$ is a direct cause of $L_2$,
the effect along green paths is not identifiable.
}
\label{mixed_time_med_fail}
\end{figure}

The proof of this theorem is given in the supplementary materials.
As an example, assume if we are interested in the effect of $A_0,A_1$ on $Y$
along only the green paths in the graph in Fig. \ref{time_med}.  The
set of nodes $V^*$ in this case is just $\{ L_1, L_2, M_1, M_2, Y \}$.  There
are three districts in the corresponding graph ${\cal G}_{V^*}$, which is just
the graph obtained from Fig. \ref{time_med} by removing $A_0,A_1$ and all
edges adjacent to these nodes.  These three districts are
$\{ M_1 \}$, $\{ M_2 \}$, and $\{ L_1, L_2, Y \}$.  It is never the case that
both a green and a blue arrow from $A_0$ or $A_1$ points to nodes in the
same district such that these nodes are ancestors of $Y$.  This means there is
no recanting district for this effect, which
in turn means the effect is expressible as a functional of
interventional densities.
We have already verified this fact in the previous section where this functional
was given as equation \ref{pot_out_example}, and which can be rephrased in
the $\text{do}(.)$ notation as follows:
\begin{equation*}
\sum_{l_1,l_2,m_1,m_2} 
	p(Y, l_1, l_2 | \text{do}(a'_0,a'_1,l_1,l_2,m_1,m_2))
	\cdot p(m_1, m_2 | \text{do}(a_0,a_1,l_1,l_2))
\end{equation*}
which can be shown equals to
\begin{equation}
\label{do_example}
\sum_{l_1,l_2,m_1,m_2} 
	p(Y, l_1, l_2 | \text{do}(a'_0,a'_1,l_1,l_2,m_1,m_2))
	\cdot p(m_1 | \text{do}(a_0,l_1)) \cdot
	\cdot p(m_2 | \text{do}(a_1,l_2))
\end{equation}
which is easily verified to be an example of equation \ref{inter_g}.

On the other hand, in either of the graphs shown in Fig.
\ref{mixed_time_med_fail}, the recanting district exists.  In Fig.
\ref{mixed_time_med_fail} (a), the district $\{ L_1, L_2, M_2, Y \}$ is
recanting, since the path $A_1 \to L_2 \to Y$ is not fully green (e.g.
we are not interested in the effect along this path), while the path
$A_1 \to M_2 \to Y$ is fully green, and $L_2$ and $M_2$ lie in the same district.
Similarly, in Fig. \ref{mixed_time_med_fail} (b), the district
$\{ M_1 \}$ is recanting, since the path $A_0 \to M_1 \to L_2 \to Y$ is
not fully green, while the path $A_0 \to M_1 \to Y$ is fully green, and
$M_1$ is its own district.

The recanting district criterion generalizes an earlier result for static
treatments and no unobserved confounding known as a \emph{recanting witness},
where the ``witness'' is a singleton node \cite{chen05ijcai}.  The term is ``recanting''
because for the purposes of one path from a particular treatment
$A_k$ the witness (or district in our case)
pretends the treatment should be active, while for the purposes of another
path from that same treatment $A_k$
the witness (or district) ``changes the story,'' and pretends the
treatment should be baseline.  Identification of path-specific effects in terms
of interventional distributions must always avoid this ``recanting'' phenomenon.
Note that even in the case of multiple longitudinal treatments, the
``recanting'' phenomenon still involves a single treatment, but spoils
the identification of the whole effect of multiple treatments.

The presence of the recanting district prevents the expression of path-specific
effects in terms of either observed or interventional data in the sense that
it is possible to construct two distinct causal models which are observationally
and interventionally indistinguishable, which are represented by the same causal
diagram, but which disagree on the value of the path-specific effect.
Furthermore, if the recanting district criterion does not exist, it is possible
to characterize cases in which the expression for the path-specific effect in
terms of interventional densities can be further expressed in terms of
observational data.

\begin{theorem}\label{obs_thm}
Let ${\cal G}$ be an ADMG with nodes $V$ representing a causal diagram with
unobserved confounders corresponding to a structural causal model.
Let $A$, $Y$ sets of nodes nodes in
${\cal G}$, and $\pi$ a subset of proper causal paths which start with a node in
$A$ and end in a node in $Y$ in ${\cal G}$.  Assume there does not exist
a recanting district for the $\pi$-specific effect of $A$ on $Y$.
Then the counterfactual representing the $\pi$-specific effect of $A$ on $Y$
is expressible in terms of the observed data $p(V)$ if and only if the total
effect $p(y | \text{do}(a))$ is expressible in terms of $p(V)$.
Moreover, the functional of $p(V)$
equal to the counterfactual is obtained from equation (\ref{inter_g}) by
replacing each interventional term in (\ref{inter_g})
by a functional of the observed data identifying that term
given by Tian's identification algorithm
\cite{tian02general},\cite{shpitser06id,shpitser07hierarchy}.
\end{theorem}

General theory of identification of causal effects
states that if
$p(y | \text{do}(a))$ is identifiable in terms of observed data, then it can
be expressed as the functional very similar to that in equation
(\ref{inter_g}), except all
variables in $A$ are assigned ``active values'' $a$.
If $p(y | \text{do}(a))$ is identifiable from observed data, then each
of the interventional terms in the functional is expressible in terms of
observed data.  This theorem simply states that to
obtain our path-specific effect all we have to do is obtain the functional
of the observed data expressing $p(y | \text{do}(a))$, and
replace the appropriate ``active values'' $a$ by ``baseline values'' $a'$
in those terms of the functional which correspond to districts containing
children of treatment $A$ via blue arrows.
For example, it can be shown that the total effect $p(y | do(a_0,a_1))$ in
Fig. \ref{time_med} is equal to
\begin{equation} \label{observational_g_}
\sum_{l_1,l_2,m_1,m_2}
	p(Y | a_0,a_1,l_1,l_2,m_1,m_2) \cdot
	p(m_2 | l_2, a_1, m_1, a_0) \cdot
	p(l_2 | a_0, a_1, l_1) \cdot
	p(m_1 | l_1, a_0) \cdot
	p(l_1 | a_0)
\end{equation}

Replacing $a_0$ and $a_1$ by $a'_0$ and $a'_1$ in expressions for
$Y$, $L_1$ and $L_2$ (which are the only terms where the node before
the conditioning bar is a child of $A_0$ or $A_1$ along blue arrows)
yields precisely equation (\ref{observational_g}) which is the function
of the observed data equal to the path-specific effect.

\begin{corollary}[generalized mediation formula for path-specific effects]
\label{med_form}
Let ${\cal G}$ be an ADMG with nodes $V$
representing a causal diagram with unobserved
confounders corresponding to a structural causal model.
Let $A$ be a set of nodes, $Y$ a single node in
${\cal G}$, and $\pi$ a subset of proper causal paths which start with a node in
$A$ and end in $Y$ in ${\cal G}$.  Assume there does not exist
a recanting district for the $\pi$-specific effect of $A$ on $Y$.
Assume $p(y | \text{do}(a))$ is expressible as a functional
$f_{\text{do}(a)}(p(V))$ of the observed data, and the path-specific
effect is equal to the functional $f_{\pi}(p(V))$ obtained in Theorem
\ref{obs_thm}.  Then the path-specific
effect along the set of paths $\pi$ on the mean difference scale for
active value $a$ and baseline value $a'$ is equal to
\[
\text{Effect along paths in $\pi$} =
E[Y]_{f_{\pi}(p(V))} - E[Y]_{f_{\text{do}(a')}(p(V))}
\]
while the path-specific effect along all paths not in $\pi$ on the mean
difference scale is equal to
\[
\text{Effect along paths not in $\pi$} =
E[Y]_{f_{\text{do}(a)}(p(V))} - E[Y]_{f_{\pi}(p(V))}
\]
\end{corollary}

An example of the generalized mediation formula applied to the longitudinal
mediation setting shown in Fig. \ref{time_med} is shown in equation
(\ref{dir_diff_gen}) for the direct effect, and equation (\ref{indir_diff_gen})
for the indirect effect.
The proof of these assertions is also given in the supplementary materials.

\section{Conclusions}

In this paper, we have shown that existing methods for mediation analysis
in epidemiology and psychology literature
based on the product and difference methods \cite{judd81process},\cite{baron86mm}
and linear
regression models suffer from problems in the presence of interactions,
non-linearities, binary outcomes, unobserved confounders, and other
modeling complications.  We have described a general framework developed in the
causal inference literature based on potential outcome counterfactuals,
non-parametric structural equations, and causal diagrams which recovers
the product and difference methods as a special case, but which is flexible
enough to handle multiple types of difficulties which arise in practical
mediation analysis situations.

Our paper serves two aims.  We first wish to caution against careless use
of mediation methodology based on linear regressions in situations where such
methodology is not suitable.  Such careless use may invalidate any conclusions
about mediation that are drawn.  Second, we want to show that appropriate use of
functional models and potential outcomes is a very flexible strategy for
tackling complex questions in causal inference, including mediation
questions in longitudinal settings with unobserved confounding.
We demonstrate this flexibility by developing a complete
characterization of situations when path-specific effects are expressible
as functionals of the observed data.  This result paves the way for using
statistical tools for answering general mediation questions in longitudinal
observational studies.
We argue that methods based on functional models
and potential outcomes are often a more appropriate methodology
in complex mediation setting than simpler methods based on linear
structural equations.

\bibliographystyle{apacite}
\bibliography{references}


\end{document}


\newcounter{dogsbody} 
\newcounter{catsbody} 
\newcounter{isamac} 
\setcounter{isamac}{0} 

\ifthenelse{\value{isamac}=1}{
\input BoxedEPS          %
\SetTexturesEPSFSpecial  
\HideDisplacementBoxes   %
}{}

\newcommand{\stt}{\small\tt}          
\newcommand{\ssy}[1]{\hbox{\tiny\itshape #1}}   


\newtheorem{theorem}{Theorem}
\newtheorem{corollary}[theorem]{Corollary}
\newtheorem{proposition}[theorem]{Proposition}
\newtheorem{lemma}[theorem]{Lemma}
\newtheorem{definition}[theorem]{Definition}
\newtheorem{remark}[theorem]{Remark}
\newtheorem{algorithm}[theorem]{Algorithm}
\newtheorem{conjecture}[theorem]{Conjecture}

\newtheorem{example}[theorem]{Example}
\newtheorem{question}[theorem]{Question}
\newtheorem{problem}[theorem]{Problem}

\newenvironment{figt}{\begin{figure}}

\newenvironment{property}[1]{\vskip10pt
{\noindent \bf #1} 
\vskip3pt\par \it      
\noindent}{\vskip10pt}      

\newcommand{\ds}{\hbox{d-sep}\-ar\-a\-tion}
\newcommand{\dsed}{\hbox{d-sep}\-ar\-at\-ed}
\newcommand{\dses}{\hbox{d-sep}\-ar\-at\-es}
\newcommand{\dc}{\hbox{d-con}\-nect\-ion}
\newcommand{\dct}{\hbox{d-con}\-nect}
\newcommand{\dcing}{\hbox{d-con}\-nect\-ing}
\newcommand{\dcs}{\hbox{d-con}\-nect\-s}
\newcommand{\dced}{\hbox{d-con}\-nect\-ed}
\newcommand{\ccon}{\hbox{co-con}\-nect\-ion}


\newcommand{\an}{{\rm an}}       
\newcommand{\ant}{{\rm ant}}     
\newcommand{\de}{{\rm de}}       
\newcommand{\spo}{\mathrm{sp}}      
\newcommand{\nsp}{\mathrm{nsp}}      
\newcommand{\ch}{{\rm ch}}       
\newcommand{\pa}{{\rm pa}}       
\newcommand{\nei}{{\rm ne}}      
\newcommand{\un}{{\rm un}}       
\newcommand{\mb}{{\rm mb}}       
\newcommand{\pre}{{\rm pre}}     

\newcommand{\path}[1]{{\ensuremath{\boldsymbol{#1}}}} 

\newcommand{\mar}[3]{\ensuremath{{\left.#1\right[}^{\scriptscriptstyle#2}
_{\scriptscriptstyle#3}}} 

\newcommand{\indm}[2]{\ensuremath{{\mathfrak I}_{\kern-1pt\scriptstyle#1}({\mathcal
#2})}} 

\newcommand{\trip}[3]{\ensuremath{\langle #1,
#2\mid #3\rangle}} 

\newcommand{\Trip}[3]{\ensuremath{\bigl\langle #1,
#2 \bigr\vert\bigl. #3\bigr\rangle}} 

\newcommand{\bigtrip}[3]{\ensuremath{\Bigl\langle #1,
#2 \Bigr\vert\Bigl. #3\Bigr\rangle}} 

\newcommand{\ord}[1]{\ensuremath{\left< #1\right>}} 

\newcommand{\ind}{\mbox{$\perp \kern-5.5pt \perp$}}
\newcommand{\nind}{\mbox{$\not\hspace{-4pt}\ind$}}

\newcommand{\mapright}[1]{\smash{
                        \mathop{\longrightarrow}\limits^{#1}}}
\newcommand{\mapdownr}[1]{\Big\downarrow 
        \rlap{$\vcenter{\hbox{$\scriptstyle#1$}}$}}
\newcommand{\mapdownl}[1]{ 
        \llap{$\vcenter{\hbox{$\scriptstyle#1$}}$}\Big\downarrow}

\newcommand{\nmodels}{\mbox{$\;\not\hspace{-2pt}\models$}}


\newcommand{\uned}{\hbox{\kern3pt\raise2.5pt\vbox{\hrule
width9pt height 0.3pt}\kern3pt}}

\newcommand{\dashed}{\hbox{\kern3.05pt\raise2.5pt\vbox{\hrule
width1.7pt height 0.3pt}\kern1.8pt\raise2.5pt\vbox{\hrule
width1.7pt height 0.3pt}\kern1.8pt\raise2.5pt\vbox{\hrule
width1.7pt height 0.3pt}\kern1.8pt\raise2.5pt\vbox{\hrule
width1.7pt height 0.3pt}\kern3.05pt}}


\newcommand{\head}{\ensuremath{\succ}}
\newcommand{\tail}{\ensuremath{-}}

\newcommand{\pedg}[2]{\ensuremath{{\kern0.5pt
\scriptstyle{\ifthenelse{\equal{\head}{#1}}{\lhead\kern0.5pt}{#1\kern0.5pt}}\joinrel\relbar
\negthinspace\relbar\joinrel{\kern0.5pt #2}\kern0.5pt}}}

\newcommand{\pdots}{\hbox{\kern2.5pt\raise1.5pt\hbox{\ensuremath{\ldots}}\ke
rn2.5pt}}  

\newcommand{\RRR}{\mathbb{R}}
\newcommand{\NNN}{\mathbb{N}}
\newcommand{\ND}{\mathcal{N}}
\newcommand{\tr}{\mathrm{tr}}
\newcommand{\bd}{\mathrm{bd}}
\newcommand{\Bd}{\mathrm{Bd}}

\newcommand{\bi}{\leftrightarrow}
\newcommand{\arr}{\mathrm{arr}}

\maketitle

\appendix

\section{Overview}

Here we give a complete theory underlying our reported results.
We proceed as follows.  We first describe some background terminology on
graphs, potential outcomes, and causal inference.
We then give a formal definition of path-specific
effects for multiple treatments, and multiple outcomes in causal diagrams where
some nodes are possibly unobserved.  This is one of the most general settings
for mediation problems.  In particular, it subsumes longitudinal mediation
settings with unobserved confounding.  We do not, however, handle issues that
arise due to selection bias, measurement error, and model misspecification bias.
These difficulties are very general in causal inference and statistics itself,
and handling them in mediation settings is a topic of future work.
We then give a formula which expresses a path-specific
effect in terms of interventional densities if the recanting district does not
exist for the effect, in accordance with Theorem \ref{full-result}.  We then
prove that if such a district does not exist, then the formula is indeed
correct (the ``if'', or soundness part of Theorem \ref{full-result}).
We then show that if a recanting
district does exist, then it is possible to construct two causal models which
agree on the causal diagram, and all interventional densities, but disagree on
the path-specific effect of interest.  This implies the path-specific effect
cannot be expressed as a functional of interventional densities for all causal
models which induce a particular diagram (the ``only if'' or completeness
part of Theorem \ref{full-result}).
Finally, we give the proofs of Theorem \ref{obs_thm}, and Corollary
\ref{med_form}.


\subsection{A Note on Notation}

We will use a somewhat different notation in this supplement from that found in
the main body of the paper, for the following reasons.
In this paper, we need to clearly distinguish mathematical
objects along the following three axes: graphical (nodes) vs statistical
(variables), singletons vs sets, and values vs variables.  In simple examples
such as those found in the main body of the paper, the singleton/set distinction
did not really arise, and it was clear from context whether graphical or
statistical concepts were discussed.  Then values can be distinguished by lower
case letters, while variables can be distinguished by upper case letters.
In the supplement, where we develop a general mathematical theory underlying the
examples, this option is no longer open to us -- the uppercase/lowercase Roman
notation is simply ambiguous when all three distinctions need to be drawn.

The paper must satisfy conflicting notation goals: we want to both make the
main body of the paper accessible, and conform to existing notation in
mediation papers, and at the same time to make the mathematical portion of the
paper unambiguous and clear.  The change of notation satisfies these two goals,
although it makes the mapping of the main body of the paper to the supplement
somewhat harder.  To assist this mapping, we restate all main results in the
supplement with the new notation.  With that said, the notation used in this
supplement will be as follows.  Vertices in a graph
${\cal G}$ will be denoted by lower case letters: $w$.  Sets of vertices
will be denoted by upper case letters: $W$.  A random variable corresponding
to a vertex $w$ will be denoted $X_w$, or sometimes subscripted:
$X_1, \ldots X_n$.  A set of random variables corresponding to a set of
vertices $W$ will be denoted by $X_W$.  A value assignment to a random variable
$X_w$ will be denoted by $x_w$.  A value assignment to a set of random variables
$X_W$ will be denoted by $x_W$.

We will consider sets of random variables $X_{v_1}, \ldots, X_{v_n}$,
and joint
probability distributions $p(x_{v_1}, \ldots, x_{v_n})$ over these sets.

\subsection{Graphs, Counterfactuals, and Structural Models of
	\\Causality}

A directed graph is a graph containing vertices (or nodes) and
directed arrows connecting pairs of vertices.  If vertices $w,y$ in a graph
${\cal G}$ are connected by a directed edge $w \to y$, we say $w$ is a parent
of $y$ and $y$ is a child of $w$.  A sequence of nodes such that every $k$th
and $k+1$th node in the sequence are connected by an edge, and no node occurs
more than once in a sequence is called a \emph{path}.  If vertices $w,y$ are
connected by a path of the form $w \to \ldots \to y$, then we say $w$ is an
ancestor of $y$ and $y$ is a descendant of $w$.  A directed graph is acyclic if
no node is its own ancestor.  We abbreviate directed acyclic graphs as DAGs.
For a given node $w$ in ${\cal G}$, we denote
its sets of parents, children, ancestors and descendants
as $Pa_{{\cal G}}(w)$, $Ch_{{\cal G}}(w)$, $An_{{\cal G}}(w)$,
$De_{{\cal G}}(w)$, respectively.  The ``genealogic relations'' on sets of
vertices are defined by taking unions, for instance for a set $W$,
$An_{{\cal G}}(W) = \bigcup_{w_i \in W} An_{{\cal G}}(w_i)$.

A Bayesian network is a DAG ${\cal G}$ which contains $n$ nodes,
$\{ v_1, \ldots v_n \} = V$, and a set of random variables $X_{v_1}, \ldots
X_{v_n}$, one for each vertex in ${\cal G}$, forming a joint probability
distribution $p(x_{v_1}, \ldots, x_{v_n})$ with a certain property linking the
distribution and the graph.  This property is called the
Markov factorization property:

$$p(x_{v_1}, \ldots, x_{v_n}) = \prod_{i = 1}^n p(x_{v_i} \mid
x_{Pa_{\cal G}(v_i)})$$

This factorization is equivalent to the local Markov property which states
that each $X_{v_i}$ is independent of $X_{V \setminus (De_{\cal G}(v_i) \cup Pa_{\cal G}(v_i))}$
conditional on $X_{Pa_{\cal G}(v_i)}$, and in turn equivalent to the global Markov
property defined by d-separation \cite{pearl88probabilistic}, which states,
for any disjoint sets of vertices $W,Y,Z$ in ${\cal G}$,
that if all paths of a certain type from nodes in $W$ to nodes in $Y$ are ``blocked'' by
nodes in $Z$ in ${\cal G}$, then $X_W$ is independent of $X_Y$ given $X_Z$ in
$p(x_{v_1}, \dots, x_{v_n})$.  A Bayesian network is thus a statistical model in that its
Markov properties define a set of probability distributions.

A statistical graph model can further be considered a causal model if we can
meaningfully talk about interventions on variables.  An intervention on $A$,
denoted by $\text{do}(X_A = x_A)$ (which we will shorten to $\text{do}(x_A)$),
is an operation that sets the variables $X_A$ to values $x_A$,
regardless of the usual behavior of $X_A$ given by the observable joint
distribution $p$.  Effects of such interventions on other variables in
the system will represent causal effects.  A statistical model defined by
a DAG ${\cal G}$ is causal if for every such $\text{do}(x_A)$,

$$
p(x_{V \setminus A} \mid \text{do}(x_A)) =
\prod_{v_i \not\in A} p(x_{v_i} \mid x_{Pa_{\cal G}(v_i)})
$$

Informally, this formula asserts that whenever we intervene on a set of
variables $X_A$, we remove from the Markov factorization all terms
$p(x_a \mid x_{Pa_{{\cal G}}(a)})$, for all $a \in A$.
This is known as the truncation formula
\cite{spirtes93causation},\cite{pearl00causality}, or
the g-formula \cite{robins86new}.  This formula implies, in particular,
that for any
$X_a$, $p(x_a \mid do(x_{Pa_{\cal G}(a)})) = p(x_a \mid x_{Pa_{\cal G}(a)})$.

The intuition for the g-formula is that in a causal model the parents of every
variable are that variable's \emph{direct causes}.  These direct causes
determine with what probability a variable assumes its values in the model.
By intervening on a variable,
we force it to attain a particular value, independently of the usual influence
of direct causes.  For this reason, we ``drop out'' the term which links the
direct causes and the variable from the Markov factorization.

We denote a random variable $X_y$ in a causal model
after $do(x_A)$ has been performed by the notation $X_y(x_A)$.  Such a variable
is called a potential outcome, or a counterfactual variable.
An assumption very commonly made in causal inference is the so called
consistency assumption, which states that if we observed variables $X_A$
attain a value $x_A$, then for any $X_Y$, the variable sets $X_Y$ and
$X_Y(x_A)$ are the same.  This assumption is crucial
in that it allows us to link outcomes under hypothetical interventions with
outcomes seen in observational studies where no interventions were in fact
performed.

In order to express natural direct and indirect effects in terms of
interventional distributions, it was necessary to assume independence of
counterfactual variables which lie in different hypothetical worlds.
This is not a testable assumption, since no possible experiment we could
perform can falsify it.  Nevertheless, there is one type of causal model that
implies such assumptions in a plausible way.  This causal model is
the so called \emph{non-parametric structural equation model}
\cite{pearl00causality} (sometimes abbreviated as an NPSEM, although
from this point on we will refer to this model simply
as a functional model.)  This model is a graphical model which consists of
a distribution $P(x_{v_1}, \ldots, x_{v_n})$ that factorizes according to a
DAG ${\cal G}$ such that every intervention can be expressed in terms of
the g-formula, and the consistency assumption is true for every
counterfactual.  In addition, we assume that every observable variable
$X_{v_i}$ is causally determined from its direct causes
$X_{Pa_{\cal G}(v_i)}$ (plus possibly a single unobserved cause only
of $X_{v_i}$ and no other variable) via some unknown function
or causal mechanism.  A functional model for the simple mediation setting
is shown in equations (\ref{y_model_general}), (\ref{m_model_general}), and
(\ref{a_model_general}).

Because of these functions, this functional model can be viewed as a kind of
``stochastic circuit'' with variables representing wires, and functions
representing logic gates that determine the voltage at a particular wire
in terms of other wires in the circuit, with a few specific wires allowed
to be randomized.  On the one hand, the model may seem quite reasonable, since
many data generating processes in nature can be naturally thought of as such
circuits (at least on the level of Newtonian physics).  On the other hand,
the model implies that for every $X_w,X_y$, it is the case that
$X_w(x_{Pa_{\cal G}(w)})$ is independent of $X_y(x_{Pa_{\cal G}(y)})$, for any
value assignments $x_{Pa_{\cal G}(w)}$,$x_{Pa_{\cal G}(y)}$, even if
$Pa_{\cal G}(w)$ and $Pa_{\cal G}(y)$ have nodes in common, and these nodes are
set to conflicting values by $x_{Pa_{\cal G}(w)}$ and $x_{Pa_{\cal G}(y)}$.
For this reason, these kinds of functional models are sometimes considered too
strong \cite{robins10alternative}.
In the remainder of the supplement, we will assume our graphs represent
functional models,
with a warning that without assuming such strong models directly, or at least
certain cross-world independences such models imply, none of the identification
results presented in this paper are valid.

\subsection{A General Definition of Path-Specific Effects}

We define path-specific effects in the general case of multiple
treatments $\{ a_1, \ldots a_k \} = A$ and multiple outcomes
$\{ y_1, \ldots y_m \} = Y$.
We will consider two sets of values for $A$, the treatment values
$\{ x_{a_1}, \ldots x_{a_k} \} = x_A$ and reference values
$\{ x^*_{a_1}, \ldots x^*_{a_k} \} = x^*_A$.
We will be interested in effect of the set $X_A$ on the set $X_Y$ along a set
$\pi$ of directed paths from nodes in $A$ to nodes in $Y$.  As described in
the main body of the paper, our
convention will be that if an arrow is a part of some path in $\pi$, it is
shown in green, otherwise it is shown in blue.

We now show how to translate path-specific effects along a given set of paths
into counterfactual form.  The resulting counterfactual, denoted by
$X_Y(\pi(x_A), x^*_A)$, will be defined inductively.

Let $V^* = An(Y)_{{\cal G}_{\not A}}$, where ${\cal G}_{\not A}$ is a subgraph
of ${\cal G}$ containing all vertices other than $A$, and all edges between
these vertices which occur in ${\cal G}$.
Fix a node $s$ in $V^*$, and let $Pa_s = Pa_{{\cal G}}(s)$.  Let $B$ be the set
of nodes $t \in Pa_s$ such that the arrow $t \to s$ is green.  For each such
$t$, let $X_t(\pi(x_A), x^*_A)$ be the inductively defined path-specific effect
of $X_A$ on $X_t$ along $\pi$.  Then the path-specific effect of $X_A$ on $X_s$
along $\pi$ is defined as

\begin{equation*}
X_s(X_{B\setminus A}(\pi(x_A), x^*_A),
X_{Pa_s \setminus (B \cup A)}(x^*_A), x_{A \cap B}, x^*_{(A \cap Pa_s)
	\setminus B})
\end{equation*}

The path-specific effect $X_Y(\pi(x_A), x^*_A)$ is defined to be the joint
distribution over
$X_{y_1}(\pi(x_A), x^*_A), \ldots, X_{y_m}(\pi(x_A), x^*_A)$.

\subsection{A Restatement of the Main Results}

In this section, we restate the main definition of the recanting
district and the three main results of the paper in terms of
the notation used in this supplement.

\begin{definition}[recanting district]
Let ${\cal G}$ be an ADMG, $A$, $Y$ sets of nodes in ${\cal G}$, and
$\pi$ a subset of proper causal paths which start with a node in $A$ and end in
a node in $Y$ in ${\cal G}$.
Let $V^*$ be the set of nodes not in $A$ which are ancestral of $Y$
via a directed path which does not intersect $A$.  Then a district $D$ in
an ADMG ${\cal G}_{V^*}$ is called
a \emph{recanting district} for the $\pi$-specific effect of $A$ on $Y$
if there exist nodes $z_i, z_j \in D$ (possibly $z_i = z_j$), $a_i \in A$,
and $y_i,y_j \in Y$ (possibly $y_i = y_j$) such that there is
a proper causal path $a_i \to z_i \to \ldots \to y_i$ in $\pi$, and a proper
causal path $a_i \to z_j \to \ldots \to y_j$ not in $\pi$.
\end{definition}

\begin{thma}{\ref{full-result}
\textbf{(recanting district criterion)}
}
Let ${\cal G}$ be an ADMG representing a causal diagram with unobserved
confounders corresponding to a structural causal model.
Let $A$, $Y$ sets of nodes nodes in ${\cal G}$, and
$\pi$ a subset of proper causal paths which start with a node in $A$ and end
in a node in $Y$ in ${\cal G}$.
Then the $\pi$-specific effect of $A$ on $Y$ is expressible in terms of
interventional densities if and only if there does
not exists a recanting district for this effect.
\end{thma}

\begin{thma}{\ref{obs_thm}}
Let ${\cal G}$ be an ADMG with nodes $V$ representing a causal diagram with
unobserved confounders corresponding to a structural causal model.
Let $A$, $Y$ sets of nodes nodes in
${\cal G}$, and $\pi$ a subset of proper causal paths which start with a node in
$A$ and end in a node in $Y$ in ${\cal G}$.  Assume there does not exist
a recanting district for the $\pi$-specific effect of $A$ on $Y$.
Then the counterfactual representing the $\pi$-specific effect of $A$ on $Y$
is expressible in terms of the observed data $p(x_V)$ if and only if the total
effect $p(x_Y | \text{do}(x_A))$ is expressible in terms of $p(x_V)$.
Moreover, the functional of $p(x_V)$
equal to the counterfactual is obtained from equation (\ref{inter_g}) by
replacing each interventional term in (\ref{inter_g})
by functional of the observed data identifying that term
given by Tian's identification algorithm
\cite{tian02general},\cite{shpitser06id},
\cite{shpitser07hierarchy}.
\end{thma}

\begin{cora}{\ref{med_form}
\textbf{(generalized mediation formula for path-specific effects)}
}
Let ${\cal G}$ be an ADMG with nodes $V$ representing a causal diagram with
unobserved confounders corresponding to a structural causal model.
Let $A$ be a set of nodes, $y$ a single node in
${\cal G}$, and $\pi$ a subset of proper causal paths which start with a node in
$A$ and end in $y$ in ${\cal G}$.  Assume there does not exist
a recanting district for the $\pi$-specific effect of $A$ on $y$.
Assume $p(x_y | \text{do}(x_A)) = f_{\text{do}(x_A)}(p(x_V))$,
and the path-specific
effect is equal to the functional $f_{\pi}(p(x_V))$ obtained in Theorem
\ref{obs_thm}.  Then the path-specific
effect along the set of paths $\pi$ on the mean difference scale for
active values $x_A$ and baseline values $x^*_A$ is equal to
\[
\text{Effect along paths in $\pi$} =
E[X_y]_{f_{\pi}(p(x_V))} - E[X_y]_{f_{\text{do}(x^*_A)}(p(x_V))}
\]
while the path-specific effect along all paths not in $\pi$ on the mean
difference scale is equal to
\[
\text{Effect along paths not in $\pi$} =
E[X_y]_{f_{\text{do}(x_A)}(p(x_V))} - E[X_y]_{f_{\pi}(p(x_V))}
\]
\end{cora}
\subsection{The Soundness Proof}

We first show soundness, namely that if a recanting district does not exist,
then the path-specific counterfactual is identifiable from interventional
distributions.  We also give an expression for this counterfactual in terms of
interventional densities.

Let $\{ v_1, \ldots v_m \} = V^* = An(Y)_{{\cal G}_{\not A}}$.
Consider the counterfactual joint probability
$p(X_{y_1}(\pi(x_A), x^*_A) = x_{y_1}, \ldots,
X_{y_m}(\pi(x_A), x^*_A) = x_{y_m})$, representing the path-specific
effect of $A$ on $Y$ along paths in $\pi$.

``Unrolling'' this counterfactual, we get the following formula:

\begin{equation}
\sum_{x_{V^* \setminus Y}} p(X_{v_1}(x_{Pa_{{\cal G}_{V^*}}(v_1)}) = x_{v_1},
\ldots X_{v_m}(x_{Pa_{{\cal G}_{V^*}}(v_m)}) = x_{v_m})
\label{eqn-1}
\end{equation}

where each value assignment $x_{Pa_{{\cal G}_{V^*}}(v_i)}$ is consistent
with $x_{v_1}, \ldots x_{v_m}$ and $x_{V^* \setminus Y}$, and the
values of $X_A$ given by the effect definition (that is if there is a green
arrow from $a \in A$ to $v_i$, then $x_{Pa_{{\cal G}_{V^*}}(v_i)}$ assigns
to $X_a$ the treatment value $x_a$ rather than the reference value $x^*_a$).

One of the assumptions that functional DAG models make is that absence of
a directed arrow from $a$ to $y$ implies fixing all observable parents of $X_y$
renders the resulting counterfactual $X_y(x_{Pa_{\cal G}(y)})$ independent of
any counterfactual $X_a(.)$, and that fixing $X_a$ will not change
$X_y(x_{Pa_{\cal G}(y)})$.


This in turn implies that in a marginal of a functional DAG model represented
by an ADMG ${\cal G}$, for any two counterfactuals $X_z(x_{Pa_{\cal G}(z)})$,
$X_w(x_{Pa_{\cal G}(w)})$, if there is no bidirected arrow from $z$ to $w$ in
${\cal G}$, then $p(X_z(x_{Pa_{\cal G}(z)}), X_w(x_{Pa_{\cal G}(w)})) =
p(X_z(x_{Pa_{\cal G}(z)})) \cdot p(X_w(x_{Pa_{\cal G}(w)}))$.  Further, any
graphical independence model, including models induced by functional DAGs obey
a property called compositionality.  A counterfactual version of
compositionality states that for any sets of counterfactual variables
$(X_A(x_{S_A}) \ci X_Y(x_{S_Y})
\mid X_Z(x_{S_Z}))$ and
$(X_W(x_{S_W}) \ci X_Y(x_{S_Y}) \mid X_Z(x_{S_Z}))$ hold, then
$(X_A(x_{S_A}) \cup X_W(x_{S_W}) \ci X_Y(x_{S_Y}) \mid X_Z(x_{S_Z}))$ also
holds.

These properties imply the that formula \ref{eqn-1} is equivalent to the
following formula

\begin{equation}
\sum_{x_{V^* \setminus Y}} \prod_{v_1, \ldots, v_k \in D \in {\cal D}({\cal G}_{V^*})}
p(X_{v_1}(x_{Pa_{{\cal G}_{V^*}}(v_1)}) = x_{v_1}, \ldots
X_{v_k}(x_{Pa_{{\cal G}_{V^*}}(v_k)}) = x_{v_k})
\label{eqn-2}
\end{equation}

where ${\cal D}({\cal G}_{V^*})$ is a set of districts in ${\cal G}_{V^*}$.
This is a decomposition of formula \ref{eqn-1} into a set of terms, one for
each district in ${\cal G}_{V^*}$

Finally, since all interventional values $x_{Pa_{{\cal G}_{V^*}}(v_i)}$
for $X_{v_i}$ involve either assignments to $A$ or assignments to variables
which appear in the district, and moreover, the intervened on values
$x_{Pa_{{\cal G}_{V^*}}(v_i)}$ are consistent with the assigned values
$x_{v_j}$ ($j \neq i$) for variables in the district, we can use the consistency
assumption to conclude formula \ref{eqn-2} is equivalent to formula \ref{eqn-3}.

\begin{equation}
\sum_{x_{V^* \setminus Y}} \prod_{D \in {\cal D}({\cal G}_{V^*})}
p(X_D = x_D \mid \text{do}(x_{Pa_{{\cal G}_{V^*}}(D) \setminus D}))
\label{eqn-3}
\end{equation}

where $x_{Pa_{{\cal G}_{V^*}}(D) \setminus D}$ is a value assignment to
$Pa(D)_{{\cal G}_{V^*}} \setminus D$ consistent
with $x_{V^* \setminus Y}$ and assignments to $X_A$ given by the effect.


This establishes our result.

\subsection{The Completeness Proof}

Next, we show completeness, namely that if a recanting district exists, then the
path-specific effect given by a counterfactual distribution\\
$p(X_{y_1}(\pi(x_A), x^*_A), \ldots, X_{y_m}(\pi(x_A), x^*_A))$ is not identifiable.
The proof will proceed as follows.

We will first show if there exists a recanting district $D$ (for a particular
$a \in A$) then the following counterfactual probability
$\gamma_1$ is not identifiable from
$P_* = \{ p(X_{V \setminus W} \mid \text{do}(x_w)) \mid W \subseteq V \}$
in the graph ${\cal G}_{D \cup \{ a \}}$:

\begin{equation}
\gamma_1 = \sum_{x_{v_i} : v_i \in D \setminus
        rh(D)_{{\cal G}_{D \cup \{ a \}}}}
                p(X_{v_1}(x_{Pa_{{\cal G}_{D \cup \{ a \}}}(v_1)}) = x_{v_1},
                \ldots
                X_{v_k}(x_{Pa_{{\cal G}_{D \cup \{ a \}}}(v_k)} = x_{v_k}))
\label{prf-step-1}
\end{equation}

where $\{ v_1, \ldots v_k \} = D$,
$rh(D)_{{\cal G}_{D\cup \{a\}}} = \{ v_i \in D \mid Ch(D)_{{\cal G}_{D \cup \{ a \}}} \cap D = \emptyset \}$, and
$x_{Pa_{{\cal G}_{D \cup \{ a \}}}(v_i)}$ for every $v^i \in D$, is a value assignment defined as
follows.

It's an assignment of values
to $Pa_{{\cal G}_{D \cup \{ a \}}}(v_i)$ that are consistent with
$x_{v_i}$ (values being summed and assigned) for nodes in
$Pa(v_i)_{{\cal G}_{D \cup \{ a \}}} \setminus \{ a \}$.
If $a \in Pa_{{\cal G}_{D \cup \{ X_i \}}(v_i)}$,
the assignment assigns to $a$ the treatment value $x_a$ if the arrow from $a$
to $v_i$ is green, and the reference value $x^*_a$ otherwise
(note that by assumption there exists both a green arrow from $a$ to a node in
$D$, and a blue arrow from $a$ to a node in $D$).

After showing the non-identifiability of $\gamma_1$, we show the
non-identifiability of a related counterfactual $\gamma_2$, defined as follows.

Fix $Y' \subseteq Y$, such that all nodes in
$rh(D)_{{\cal G}_{D \cup \{a\}}}$ are ancestral of $Y'$ in
${\cal G}_{V^*}$, and for no subset of $Y'$ is this true.
For every node $r$ in $rh(D)_{{\cal G}_{D \cup \{a\}}}$ pick a node
$y_r \in Y'$ such
that there is a directed path $\pi_r$ from $r$ to $y_r$.  Let the set of
nodes in every such path be equal to $W^*$.
Let ${\cal G}^*$ be a subgraph of ${\cal G}_{V^*}$ containing nodes in
$D \cup W^*$.  Then we will show a counterfactual probability $\gamma_2$
is not identifiable from $P_*$ in ${\cal G}^*$, where $\gamma_2$ is defined as

\begin{equation}
\gamma_2 = \sum_{x_{v_i} : v_i \in (D \cup W^*) \setminus Y'}
p(X_{v_1}(x_{Pa_{{\cal G}^*}(v_1)}) = x_{v_1},
        \ldots, X_{v_l}(x_{Pa_{{\cal G}^*}(v_l)}) = x_{v_l})
\label{prf-step-2}
\end{equation}

where $\{ v_1, \ldots v_l \} = D \cup W^*$, and
$x_{Pa_{{\cal G}^*}(v_i)}$ is defined as before.

Having shown $\gamma_2$ is not
identifiable in ${\cal G}^*$ from $P_*$, we then have two models $M_1, M_2$
which agree on $P_*$ but disagree on $\gamma^*$.  We then note that augmenting
$M_1, M_2$ with additional variables can result in models $M'_1, M'_2$ that
induce ${\cal G}$, and such that $\gamma_2$ is a marginal distribution of
the counterfactual $\gamma$ in these models.  This will imply
$\gamma$ is not identifiable from $P_*$ in ${\cal G}$, which was what we
wanted to show.


\begin{lemma}
The counterfactual $\gamma_1$ given in equation \ref{prf-step-1} is not
identifiable from $P_*$ in ${\cal G}_{D \cup \{ a \}}$.
\label{lemma-1}
\end{lemma}
\begin{prf}
Pick two nodes in $D$, $v_1, v_2$ such that $a$ has a green arrow to $v_1$,
and a blue arrow to $v_2$.  Assume without loss of generality that $a$ only
affects those two nodes in $D$.  Assume, also without loss of generality,
that every node in $D$ has at most one child (other arrows are vacuous).

We now construct two functional models, ${\cal M}_1$ and ${\cal M}_2$, which both agree on $P_*$,
both induce ${\cal G}_{D \cup \{ a \}}$, but which disagree on $\gamma_1$
as defined.  In these models, every observable variable is binary.  Every
bidirected arc corresponds to an unobserved binary variable with two children.
In ${\cal M}_1$, for every observable node, its value is determined by the bit parity
function of its parents (both observed and unobserved).  For ${\cal M}_2$, for every
observable node, its value is determined by the bit parity function of
its parents, except the functions determining the values of $v_1, v_2$ do not
take the value of $a$ into account.  The distributions over unobserved
nodes is the same in both models, and is uniform.

We now show the two models have the desired properties.  That both models
induce ${\cal G}_{D \cup \{ a \}}$ is clear.
Next, we show ${\cal M}_1$ and ${\cal M}_2$ agree on $P_*$.

By construction, both models agree on $p(x_a)$.  We next show both models
agree on $p(x_D \mid \text{do}(x_a))$.
It's not difficult to show (following the proof of Theorem 17 in
\cite{shpitser07hierarchy}) that
$p(x_D \mid \text{do}(x_a)) = p(x_D)$ is a uniform distribution in ${\cal M}_2$
over assignments $x_D$ such that $x_{rh(D)_{{\cal G}_{D \cup \{a \}}}}$ has
even bit parity.  In fact, the same proof shows the same for
$p(x_D \mid \text{do}(a))$ in ${\cal M}_1$.  This implies that
$p(x_{D\cup \{ a \}}) = p(x_D \mid \text{do}(x_a)) p(x_a)$ is the same in
${\cal M}_1$ and ${\cal M}_2$.  Furthermore, for any partition $(Z_1, Z_2)$ of
$Z = (D \cup \{ a \})$, it is either the case
that $Z_1 \subset D$, or $p(x_{Z_1} \mid \text{do}(x_{Z_2})) =
p(x_{Z_1 \setminus \{ a \}} \mid \text{do}(x_{Z_2 \cup \{ a \}})) p(x_a)$.
In the former case we have two causal
submodels derived from ${\cal M}_1, {\cal M}_2$ which agree on functions and distributions
of unobserved variables, and which have the observed distribution
$p(x_D \mid \text{do}(x_a))$.  This implies ${\cal M}_1$ and ${\cal M}_2$ must agree on
$p(x_{Z_1} \mid \text{do}(x_{Z_2}))$.
In the latter case, the decomposition of the effect, and the previous argument
implies our conclusion.

Finally, we must show ${\cal M}_1$ and ${\cal M}_2$ disagree on $\gamma_1$.

In ${\cal M}_2$, $\gamma_1$ is a distribution over nodes in
$R = rh(D)_{{\cal G}_{D \cup \{ a \}}}$.  By assumption, the values of the
variables in set $X_R$ can be viewed as giving the bit parity of each unobserved
value, counted twice.  This implies $\gamma_1$ assigns probability $0$ to
any value assignment to $X_R$ with odd bit parity, and a uniform distribution
to even bit parity assignments.  What we must now show is that $\gamma_1$ is a
different distribution in ${\cal M}_1$.  Indeed, in ${\cal M}_1$ the values of
the variables in set $X_R$ can be viewed as giving the bit parity of each
unobserved value, counted twice, plus 1 (because $a$ has exactly one
directed path to $R$ in ${\cal G}_{D \cup \{ a \}}$ where $a$ takes on
value $x_a = 1$, and exactly one directed path to $R$ in
${\cal G}_{D \cup \{ a \}}$ where $a$ takes on value $x_a = 0$).
This implies $\gamma_1$ assigns probability $0$ to any value assignment
to $X_R$ with even bit parity, and a uniform distribution to odd bit parity
assignments.

The constructed models ${\cal M}_1, {\cal M}_2$ induce non-positive probabilities
$p(x_{D \cup \{ a \}})$.  It is not difficult to augment these models to create a pair
of new models ${\cal M}'_1, {\cal M}'_2$ such that $p(x_{D \cup \{a\}})$ in the new models is
positive, and the models agree on $P'_*$ (the set of interventional distributions in these
new models) and disagree on $\gamma_1$.

We construct ${\cal M}'_1, {\cal M}'_2$ by adding a new unobserved binary parent for every
node in $R$, with a distribution $\{ \epsilon, 1 - \epsilon \}$, where
$\epsilon$ is a very small positive real number.  Clearly, ${\cal M}'_1, {\cal M}'_2$
agree on any member of $P'_*$ involving nodes in
$(D \cup \{ a \}) \setminus R$.  Note that any member $P'_j$ of $P'_*$
involving nodes $R' \subseteq R$ in ${\cal M}'_1, {\cal M}'_2$ is a function
of some interventional distribution over parents of $R'$,
the distribution $P(x_{U_R})$ over unobserved
parents $U_R$ of $R$ added to ${\cal M}'_1, {\cal M}'_2$, the functions
determining the values of $R$ in ${\cal M}'_1, {\cal M}'_2$, and the
distribution over original unobserved nodes in ${\cal M}'_1, {\cal M}'_2$.
Since ${\cal M}'_1, {\cal M}'_2$ agree on all these objects,
they must agree on $P'_*$.

By construction, the probability of $\gamma_1$ in ${\cal M}'_2$ assigns low
but non zero probabilities to odd bit parity assignments to $X_R$, while the
probability of $\gamma_1$ in ${\cal M}'_1$ assigns low but non zero
probabilities to even parity assignments to $X_R$.  Since $\epsilon$ can be
made arbitrarily small, this implies ${\cal M}'_1, {\cal M}'_2$ disagree on
$\gamma_1$.

This concludes our proof.
\end{prf}

\begin{lemma}
The counterfactual $\gamma_2$ shown in equation \ref{prf-step-2} is not
identifiable from $P_*$ in ${\cal G}^*$.
\label{lemma-2}
\end{lemma}
\begin{prf}
Without loss of generality, assume every node in ${\cal G}^*$ has at most
one child.  Then we augment ${\cal M}'_1, {\cal M}'_2$ constructed in the proof of Lemma
\ref{lemma-1} by adding a binary node for every vertex in ${\cal G}^*$, but not
${\cal G}_{D \cup \{ a \}}$.  We let each such node obtain its value from
the bit parity of its parents in ${\cal G}^*$
(without adding unobserved parents).  Call the resulting models ${\cal M}"_1, {\cal M}"_2$.

Every node added to ${\cal M}"_1,{\cal M}"_2$ forms its own district, and for every such node
$w$, the distribution $p(x_w \mid \text{do}(x_{Pa(w)_{{\cal G}^*}}))$ is the same in
${\cal M}"_1$ and ${\cal M}"_2$ by construction.  This implies ${\cal M}"_1,{\cal M}"_2$ agree on
$P"_*$.
But by construction we also obtain that $M"_1, M"_2$ disagree on $\gamma_2$.

As before, the constructed models $M"_1, M"_2$ do not yield positive
observable distributions.  We augment our models and create a new
pair of models ${\cal M}^*_1, {\cal M}^*_2$ which induce positive observable distributions,
which agree on $P_*$ and disagree on $\gamma_2$.  To do so, we add for
every node in $W^* \setminus rh(D)_{{\cal G}^*}$ a new binary unobserved
parent with probabilities $\{ \epsilon, 1 - \epsilon \}$, where $\epsilon$
is a very small positive real number.  Since every node $w$ in
$W^* \setminus rh(D)_{{\cal G}^*}$ is its own district, by construction
${\cal M}^*_1,{\cal M}^*_2$ agree on $p(x_w \mid \text{do}(x_{Pa_{{\cal G}^*}(w)}))$, which implies
${\cal M}^*_1,{\cal M}^*_2$ agree on $P_*$.

The probability of $\gamma_2$ in ${\cal M}^*_2$ then assigns a small but
positive probability to any even bit parity assignment to $Y'$, while the
probability of $\gamma_2$ in ${\cal M}^*_1$ assigns a small but positive
probability to any odd bit parity assignment to $Y'$.  Since $\epsilon$ can be
made arbitrarily small, this implies ${\cal M}^*_1, {\cal M}^*_2$ disagree on
$\gamma_2$.

This establishes our result.
\end{prf}

\begin{lemma}
The counterfactual $\gamma$ is not
identifiable from $P_*$ in ${\cal G}$.
\label{lemma-3}
\end{lemma}
\begin{prf}
This can be easily established by augmenting models ${\cal M}^*_1, {\cal M}^*_2$ constructed
in the previous Lemma with enough extra nodes to enlarge ${\cal G}^*$ to
${\cal G}$.  These extra nodes will be fully jointly independent of each other and nodes in ${\cal G}^*$.  (That is, any edge connecting to such
nodes in ${\cal G}$ will be vacuous in our augmentation of ${\cal M}^*_1, {\cal M}^*_2$.
It's clear from this construction that the resulting augmented models
agree on $P_*$, disagree on $\gamma$, and induce a positive observable
distribution.

This establishes completeness of the criterion.
\end{prf}

\subsection{Remaining Proofs}

We first prove Theorem \ref{obs_thm}.
If the recanting district does not exist for a given path-specific effect
of $A$ on $Y$ along paths in $\pi$, then the corresponding counterfactual is
identifiable from interventional densities via equation (\ref{eqn-3}).
A general theory of identification \cite{tian02general},\cite{shpitser06id},
\cite{shpitser07hierarchy} states that if a causal effect
$p(x_Y | \text{do}(x_A))$ is expressible as a function of $p(x_V)$ in all
causal models inducing ${\cal G}$, then it can be expressed as a functional
almost identical to equation (\ref{eqn-3}) in this supplement,
except all values of $A$ in the
functional are $x_A$.  The conclusion immediately follows.
Corollary \ref{med_form} follows immediately from Theorem \ref{obs_thm}, and
definitions in equations
(\ref{path_specific_pi}) and (\ref{path_specific_not_pi}).

\bibliographystyle{plain}
\bibliography{references}